# SMOOTH BACKFITTING IN GENERALIZED ADDITIVE MODELS

By Kyusang Yu,[1] Byeong U. Park[2] and Enno Mammen[1]

*University of Mannheim, Seoul National University*
*and University of Mannheim*

Generalized additive models have been popular among statisticians and data analysts in multivariate nonparametric regression with non-Gaussian responses including binary and count data. In this paper, a new likelihood approach for fitting generalized additive models is proposed. It aims to maximize a smoothed likelihood. The additive functions are estimated by solving a system of nonlinear integral equations. An iterative algorithm based on smooth backfitting is developed from the Newton–Kantorovich theorem. Asymptotic properties of the estimator and convergence of the algorithm are discussed. It is shown that our proposal based on local linear fit achieves the same bias and variance as the oracle estimator that uses knowledge of the other components. Numerical comparison with the recently proposed two-stage estimator [*Ann. Statist.* **32** (2004) 2412–2443] is also made.

**1. Introduction.** In this paper, we consider *generalized additive models* where the conditional mean $m(\mathbf{x}) \equiv E(Y|\mathbf{X} = \mathbf{x})$ of a response $Y$ given a $d$-dimensional covariate vector $\mathbf{X} = \mathbf{x}$ is modeled via a known *link* $g$ by a sum of unknown component functions $\eta_i$:

$$(1) \qquad g(m(\mathbf{x})) = \eta_0 + \eta_1(x_1) + \cdots + \eta_d(x_d).$$

By employing a suitable link $g$, it allows wider applicability than ordinary additive models where $m(\mathbf{x}) = m_0 + m_1(x_1) + \cdots + m_d(x_d)$. For example, in the case where the conditional distribution of the response is Bernoulli, the conditional mean $m(\mathbf{x})$, which in this case, is the conditional probability, may be successfully modeled by a generalized additive model with the

Received June 2005; revised March 2007.

[1]Supported by DFG Project MA 1026/6-2 of the Deutsche Forschungsgemeinschaft.

[2]Supported by the Korea Research Foundation Grant funded by the Korean Government (MOEHRD) (KRF-2005-070-C00021).

logistic link $g(u) = \log\{u/(1-u)\}$. The model (1) inherits the structural simplicity and the easy interpretability of linear models. Furthermore, generalized additive models (and also additive models) are known to free one from the curse of dimensionality. Under the (generalized) additive models, one can construct an estimator of $m(\mathbf{x})$ that achieves the same optimal rate of convergence for general $d$ as for $d = 1$, see Stone [23, 24].

There have been a number of proposals for fitting the ordinary additive models. Friedman and Stuetzle [6] introduced a backfitting algorithm, and Buja, Hastie and Tibshirani [2] studied its properties. Opsomer and Ruppert [22] and Opsomer [21] showed that the backfitting estimator is well-defined asymptotically when the stochastic dependence between covariates is "not far" from independence. Mammen, Linton and Nielsen [15] proposed the so called *smooth backfitting* by employing the projection arguments of Mammen et al. [16]. In contrast to the ordinary backfitting, the dependence between covariates affects the convergence and stability of the algorithm only weakly. This was illustrated by very convincing simulations in Nielsen and Sperlich [20], where also surprisingly good performance of smooth backfitting was reported for very high dimensions. Furthermore, the local linear smooth backfitting estimator achieves the same bias and variance as the oracle estimator based on knowing the other components, and thus improves on the ordinary backfitting.

The local scoring backfitting (Hastie and Tibshirani [7]) is one of the most popular methods for generalized additive models (1). However, its theoretical properties are not well understood since it is only defined implicitly as the limit of a complicated iterative algorithm. Recently, there have been proposed other methods of fitting generalized additive models. Among others, Kauermann and Opsomer [9] proposed a local likelihood estimator which is a solution of a very large set of nonlinear score equations. They suggested an iterative backfitting algorithm to approximate the solution of the system. However, their theoretical developments are based on the assumption that the backfitting algorithm converges. Horowitz and Mammen [8] proposed a two-stage estimation procedure using the squared error loss with a link function; see also Linton [13]. In the context of local quasilikelihood estimation (see, e.g., Fan, Heckman and Wand [5]), this amounts to modelling the conditional variance to be a constant. Estimation by penalized B-splines in generalized additive models and in some related models was discussed in Eilers and Marx [4].

In this paper, we propose new estimation procedures for generalized additive models (1) that are based on a quasilikelihood with a general link. Using quasilikelihoods for fitting generalized linear models is well justified. Its advantages are similar to what maximum likelihood estimation has over other methods such as least squares approaches. The advantages carry over to the problem of fitting generalized additive models. For example, in the cases where the conditional distribution belongs to an exponential family, it guarantees convexity of the objective function if one uses the canonical link, and leads to an estimator which has the smallest asymptotic variance.

The proposed estimators solve a set of smoothed quasilikelihood equations. Unlike the least squares smooth backfitting of Mammen, Linton and Nielsen [15] in the ordinary additive models, it is a system of nonlinear integral equations. The approach is a natural generalization of parametric quasilikelihood estimation. The theoretical contribution of this paper is to show how the parametric asymptotic theory can be carried over to a nonparametric nonlinear model with several nonparametric components. The nonlinear backfitting integral equations for updating the estimators cannot be solved explicitly. This complicates a great deal development of a backfitting algorithm and its theory. We tackle this problem by employing a double iteration scheme which consists of *inner* and *outer* iterations. The outer loop is originated from a linear approximation of the smoothed quasilikelihood equations. Each step in the outer iteration is shown to be equivalent to a projection onto a Hilbert space equipped with a smoothed squared error norm, so that for each outer step we can devise a smooth backfitting procedure (inner iteration) whose limit defines an outer update. We note that the Hilbert space and its norm for each step of the outer iteration are also updated. We show that the convergence of the inner iteration is uniform for all outer loops. We discuss the smoothed quasilikelihood estimation for Nadaraya–Watson smoothing and for local linear fit. We present their theoretical properties. We find that our estimators achieve the optimal univariate rate for all dimensions. In particular, the local linear smoothed quasilikelihood estimator has the oracle bias as well as the oracle variance. Our numerical experiments also suggest that the new proposal has quite good mean squared error properties. As our estimators are defined through a projection onto an appropriate Hilbert space as the smooth backfitting technique in additive models, it is expected from the results of Nielsen and Sperlich [20] that they are successful for very high dimensions and for correlated covariates. The latter point will be illustrated by simulations in Section 5.

Some other related works on additive or generalized additive models include the marginal integration approaches of Linton and Nielsen [12], and Linton and Härdle [11]. The methods, however, suffer from the curse of dimensionality and fails to achieve the optimal univariate rate for general dimension unless the smoothness of the underlying component functions increases with dimension. See Lee [10] for a discussion on this. Mammen and Nielsen [17] considered a general class of nonlinear regression and discussed some estimation principles including the smooth backfitting. Mammen and Park [18] proposed several bandwidth selection methods for smooth backfitting, and Mammen and Park [19] provided a simplified version of the local linear smooth backfitting estimator in additive models.

The rest of the paper is structured as follows. In Section 2 we introduce the smoothed quasilikelihood estimation based on Nadaraya–Watson smoothing, and in Section 3 we extend it to the local linear framework. In Section 4 we present the asymptotic and algorithmic properties of the estimators. In Section 5 we provide the results of some numerical experiments including a comparison with the two-stage procedure of Horowitz and Mammen [8]. Finally, we give proofs and technical details in Section 6.

**2. Estimation with Nadaraya–Watson-type smoothers.** Let $Y$ and $\mathbf{X} = (X_1, \ldots, X_d)$ be a random variable and a random vector of dimension $d$, respectively and let $(\mathbf{X}^1, Y^1), \ldots, (\mathbf{X}^n, Y^n)$ be a random sample drawn from $(\mathbf{X}, Y)$. Assume that $\mathbf{X}$ has the density function $p(\cdot)$ and $X_j$ have marginal density functions $p_j(\cdot)$, $j = 1, \ldots, d$. We consider the following generalized additive model:

$$E(Y|\mathbf{X} = \mathbf{x}) = g^{-1}(\eta_0 + \eta_1(x_1) + \cdots + \eta_d(x_d)),$$

where $g$ is some known link function, $\mathbf{x} = (x_1, \ldots, x_d)$ are given value of the covariates, $\eta_0$ is an unknown constant and $\eta_j(\cdot)$, $j = 1, \ldots, d$, are univariate unknown smooth functions. Suppose that the conditional variance is modeled as $\mathrm{var}(Y|\mathbf{X} = \mathbf{x}) = V(m(\mathbf{x}))$ for some positive function $V$. The quasilikelihood function, which can replace the conditional log-likelihood when the latter is not available, equals $Q(m(\mathbf{x}), y)$, where $\partial Q(m, y)/\partial m = (y - m)/V(m)$. Note that the log-likelihood of an exponential family is a special case of a quasilikelihood function $Q(m(\mathbf{x}), y)$. The results presented in this paper for a quasilikelihood are thus valid for exponential family cases, also.

2.1. *The smoothed quasilikelihood.* Before introducing the smoothed quasilikelihood, we briefly go over the smooth backfitting in additive models proposed by Mammen, Linton and Nielsen [15]. For a Nadaraya–Watson type smoother, it starts with embedding the response vector $\mathbf{Y} = (Y^1, \ldots, Y^n)$ into the space of tuples of $n$ functions, $\mathcal{F} = \{(f^1, \ldots, f^n) : f^i \text{ are functions from } \mathbb{R}^d \text{ to } \mathbb{R}\}$. Let $K^0$ be a base kernel function and $K_h^0(u) = h^{-1}K^0(h^{-1}u)$. Define a boundary corrected kernel function by

$$(2) \qquad K_h(u, v) = \frac{K_h^0(u - v)}{\int_0^1 K_h^0(w - v)\, dw} I(u, v \in [0, 1]).$$

The space $\mathcal{F}$ is endowed with the (semi)norm

$$\|\mathbf{f}\|_*^2 = \int n^{-1} \sum_{i=1}^n (f^i(\mathbf{x}))^2 \prod_{j=1}^d K_{h_j}(x_j, X_j^i)\, d\mathbf{x}.$$

The tuple $\widetilde{\mathbf{m}} = (\widetilde{m}, \ldots, \widetilde{m})$, where $\widetilde{m}$ is the full dimensional local constant estimator, is then the projection of $\mathbf{Y}$ onto $\mathcal{F}_{\text{full}} = \{f \in \mathcal{F} : f^i \text{ does not depend on } i\}$.

The smooth backfitting estimator, denoted by $\widehat{m}$, in the form of $\widehat{\mathbf{m}} = (\widehat{m}, \ldots, \widehat{m})$ is defined as the further projection of the full dimensional estimator onto

$$\mathcal{F}_{\text{add}} = \{f \in \mathcal{F}_{\text{full}} : f^i(\mathbf{x}) \overset{i}{\equiv} g_1(x_1) + \cdots + g_d(x_d)$$

$$\text{for some functions } g_j : \mathbb{R} \to \mathbb{R}\}.$$

For tuples of functions $\mathbf{f} = (f, \ldots, f)$ in $\mathcal{F}_{\text{full}}$, one has $\|\mathbf{f}\|_*^2 = \int f(\mathbf{x})^2 \widehat{p}(\mathbf{x}) \, d\mathbf{x}$ where $\widehat{p}(\mathbf{x}) = n^{-1} \sum_{i=1}^n K_{\mathbf{h}}(\mathbf{x}, \mathbf{X}^i)$ and $K_{\mathbf{h}}(\mathbf{x}, \mathbf{X}^i) = \prod_{j=1}^d K_{h_j}(x_j, X_j^i)$. This means that $\widehat{m}(\mathbf{x}) = \widehat{m}_1(x_1) + \cdots + \widehat{m}_d(x_d)$ is the projection, in the space $L_2(\widehat{p})$, of $\widetilde{m}$ onto the subspace of additive functions $\{m \in L_2(\widehat{p}) : m(\mathbf{x}) = m_1(x_1) + \cdots + m_d(x_d)\}$. The smooth backfitting estimator $\widehat{m}$ can be obtained by projecting $\mathbf{Y}$ directly onto $\mathcal{F}_{\text{add}}$.

The smooth backfitting can be regarded as a minimization of an empirical version of $E(Y - f(\mathbf{X}))^2 = \int E[(Y - f(\mathbf{x}))^2 | \mathbf{X} = \mathbf{x}] p(\mathbf{x}) \, d\mathbf{x}$. To see this, we note that

$$\|\mathbf{Y} - f(\cdot)\mathbf{1}\|_*^2 = \int \left[ \frac{n^{-1} \sum_{i=1}^n (Y^i - f(\mathbf{x}))^2 K_{\mathbf{h}}(\mathbf{x}, \mathbf{X}^i)}{\widehat{p}(\mathbf{x})} \right] \widehat{p}(\mathbf{x}) \, d\mathbf{x},$$

where $\mathbf{1} = (1, \ldots, 1)$. This motivates us to consider the *expected quasilikelihood*, $E[Q(g^{-1}(\eta(\mathbf{X})), Y)]$, as an objective function in generalized additive models, where $\eta(\mathbf{x}) = \eta_0 + \eta_1(x_1) + \cdots + \eta_d(x_d)$. Our new estimator aims to maximize the expected quasilikelihood. This maximization can be interpreted as maximizing the quasilikelihood for all possible future observations on average.

We estimate $E[Q(g^{-1}(\eta(\mathbf{X})), Y) | \mathbf{X} = \mathbf{x}]$ by

$$\widehat{Q}_c(\mathbf{x}, \eta) = \widehat{p}(\mathbf{x})^{-1} n^{-1} \sum_{i=1}^n Q(g^{-1}(\eta(\mathbf{x})), Y^i) K_{\mathbf{h}}(\mathbf{x}, \mathbf{X}^i).$$

We use nonnegative boundary corrected kernels [see (2)], so that

$$\int K_h(u, v) \, du = 1 \quad \text{and} \quad \int K_{\mathbf{h}}(\mathbf{x}, \mathbf{X}^i) \, d\mathbf{x}_{-j} = K_{h_j}(x_j, X_j^i)$$

for $j = 1, \ldots, d$. Here and throughout the paper, $\mathbf{x}_{-j}$ denotes the vector $\mathbf{x}$ with the $j$th component $x_j$ being deleted. With a general link $g$, we define a *smoothed quasilikelihood* $SQ(\eta)$, as an estimator of the expected quasilikelihood $EQ(g^{-1}(\eta(\mathbf{X})), Y) = \int E[Q(g^{-1}(\eta(\mathbf{X})), Y) | \mathbf{X} = \mathbf{x}] p(\mathbf{x}) \, d\mathbf{x}$, by

$$(3) \qquad \begin{aligned} SQ(\eta) &= \int \widehat{Q}_c(\mathbf{x}, \eta) \widehat{p}(\mathbf{x}) \, d\mathbf{x} \\ &= \int n^{-1} \sum_{i=1}^n Q(g^{-1}(\eta(\mathbf{x})), Y^i) K_{\mathbf{h}}(\mathbf{x}, \mathbf{X}^i) \, d\mathbf{x}. \end{aligned}$$

The $L_2(\widehat{p})$ error in Mammen, Linton and Nielsen [15] is a special case of the smoothed quasilikelihood given at (3) with $Q(m, y) = -(y - m)^2/2$ and the identity link, $g(m) \equiv m$.

2.2. *Backfitting equations.* Suppose that the quasilikelihood $Q(g^{-1}(\eta), y)$ is strictly concave as a function of $\eta$ for each $y$. Since it satisfies the (conditional) Bartlett identities, $E(Q(g^{-1}(\eta(\mathbf{x})), Y)|\mathbf{X} = \mathbf{x})$ is not monotone in $\eta(\mathbf{x})$ for every $\mathbf{x}$ and thus has a unique maximizer. This implies that $SQ$ defined at (3) has a unique maximizer with probability tending to one. Let $\widehat{\eta}$ be a maximizer of $SQ(\eta)$ given at (3) over all additive functions. Then, the estimator $\widehat{\eta} = \widehat{\eta}_0 + \widehat{\eta}_1(x_1) + \cdots + \widehat{\eta}_d(x_d)$ satisfies

$$dSQ(\eta; g) = 0$$
(4)
$$\text{for all additive functions } g(\mathbf{x}) = g_0 + g_1(x_1) + \cdots + g_d(x_d),$$

where $dSQ(\eta; g)$ is the Fréchet differential of the functional $SQ$ at $\eta$ with increment $g$, see Section 7.4 in Luenberger [14]. The equation (4) is equivalent to the following set of equations:

$$\int n^{-1} \sum_{i=1}^{n} \left[ \frac{Y^i - g^{-1}(\eta(\mathbf{x}))}{V(g^{-1}(\eta(\mathbf{x})))g'(g^{-1}(\eta(\mathbf{x})))} \right] K_{\mathbf{h}}(\mathbf{x}, \mathbf{X}^i) \, d\mathbf{x} = 0,$$

$$\int n^{-1} \sum_{i=1}^{n} \left[ \frac{Y^i - g^{-1}(\eta(\mathbf{x}))}{V(g^{-1}(\eta(\mathbf{x})))g'(g^{-1}(\eta(\mathbf{x})))} \right] K_{\mathbf{h}}(\mathbf{x}, \mathbf{X}^i) \, d\mathbf{x}_{-j} = 0, \qquad j = 1, \ldots, d.$$

Let $\boldsymbol{\eta}(\mathbf{x})$ denote a tuple of functions $(\eta_0, \eta_1(x_1), \ldots, \eta_d(x_d))$. This should not be confused with $\eta(\mathbf{x}) = \eta_0 + \eta_1(x_1) + \cdots + \eta_d(x_d)$. Define

$$\widehat{F}_0\boldsymbol{\eta} = \int \left[ \frac{\widetilde{m}(\mathbf{x}) - g^{-1}(\eta(\mathbf{x}))}{V(g^{-1}(\eta(\mathbf{x})))g'(g^{-1}(\eta(\mathbf{x})))} \right] \widehat{p}(\mathbf{x}) \, d\mathbf{x},$$

$$(\widehat{F}_j\boldsymbol{\eta})(x_j) = \int \left[ \frac{\widetilde{m}(\mathbf{x}) - g^{-1}(\eta(\mathbf{x}))}{V(g^{-1}(\eta(\mathbf{x})))g'(g^{-1}(\eta(\mathbf{x})))} \right] \widehat{p}(\mathbf{x}) \, d\mathbf{x}_{-j}, \qquad j = 1, \ldots, d,$$

$$(\widehat{\mathbf{F}}\boldsymbol{\eta})(\mathbf{x}) = (\widehat{F}_0\boldsymbol{\eta}, (\widehat{F}_1\boldsymbol{\eta})(x_1), \ldots, (\widehat{F}_d\boldsymbol{\eta})(x_d))^T,$$

where $\widetilde{m}(\mathbf{x}) = \widehat{p}(\mathbf{x})^{-1}n^{-1}\sum_{i=1}^{n} Y^i K_{\mathbf{h}}(\mathbf{x}, \mathbf{X}^i)$ is the full dimensional local constant estimator. Then, $\widehat{\eta}(\mathbf{x})$ can be obtained by solving $\widehat{\mathbf{F}}\boldsymbol{\eta} = \mathbf{0}$. The estimator $\widehat{\eta}$ aims at the true $\eta^* = g(m(\cdot))$ which maximizes $\int E[Q(g^{-1}(\eta(\mathbf{x})), Y)|\mathbf{X} = \mathbf{x}]p(\mathbf{x}) \, d\mathbf{x}$, over all additive functions $\eta$.

We need to put some norming constraints on component functions for a unique identification of $\widehat{\eta}_j$ that give $\widehat{\eta}(\mathbf{x}) = \widehat{\eta}_0 + \widehat{\eta}_1(x_1) + \cdots + \widehat{\eta}_d(x_d)$. This should be done also for the component functions comprising $\eta^*$. Let $q_j(u, y)$ be the $j$th derivative of $Q(g^{-1}(u), y)$ with respect to $u$. Define for a function

$\mu$ on $\mathbb{R}^d$,

$$w^\mu(\mathbf{x}) = -q_2(\mu(\mathbf{x}), m(\mathbf{x}))p(\mathbf{x}),$$

$$\widehat{w}^\mu(\mathbf{x}) = -n^{-1}\sum_{i=1}^n q_2(\mu(\mathbf{x}), Y^i)K_\mathbf{h}(\mathbf{x}, \mathbf{X}^i).$$

We note that $w^{\eta^*}(\mathbf{x}) = g'(m(\mathbf{x}))^{-2}V(m(\mathbf{x}))^{-1}p(\mathbf{x})$ since $m(\mathbf{x}) = g^{-1}(\eta^*(\mathbf{x}))$. The function $w^\mu$ is positive for all $\mu$ if we assume $q_2(u, y) < 0$ for $u \in \mathbb{R}$ and $y$ in the range of the response. The assumption $q_2(u, y) < 0$, which is also made in Fan, Heckman and Wand [5], guarantees strict concavity of the quasilikelihood.

Let $\boldsymbol{\eta}^* \equiv (\eta_0^*, \eta_1^*, \ldots, \eta_d^*)$ maximize

$$(5) \qquad \int E[Q(g^{-1}(\eta(\mathbf{x})), Y)|\mathbf{X} = \mathbf{x}]p(\mathbf{x})\,d\mathbf{x}$$

$$\text{subject to} \quad \int \eta_j(x_j)w^\eta(\mathbf{x})\,d\mathbf{x} = 0, \qquad 1 \le j \le d.$$

If $Q(m, y) = -(y - m)^2/2$, then the norming constraints, $\int \eta_j(x_j)w^\eta(\mathbf{x})\,d\mathbf{x} = 0, 1 \le j \le d$, reduces to the usual centering condition that every component function has mean zero. We define the *maximum smoothed quasilikelihood estimator* $\widehat{\boldsymbol{\eta}}(\mathbf{x}) = (\widehat{\eta}_0, \widehat{\eta}_1(x_1), \ldots, \widehat{\eta}_d(x_d))$ to be the solution of

$$(6) \qquad \widehat{\mathbf{F}}\boldsymbol{\eta} = \mathbf{0} \quad \text{subject to} \quad \int \eta_j(x_j)\widehat{w}^\eta(\mathbf{x})\,d\mathbf{x} = 0, \qquad 1 \le j \le d.$$

2.3. *Iterative algorithms.* The major hurdle in solving $\widehat{\mathbf{F}}\boldsymbol{\eta} = \mathbf{0}$ is that it is a nonlinear system of equations, as opposed to the smooth backfitting in additive models. The approach we take to resolve this difficulty is to employ a double iteration scheme which consists of *inner* and *outer* iterations. To describe the procedure, we introduce several relevant function spaces. For a nonnegative function $w$ defined on $\mathbb{R}^d$, let $w_j$ and $w_{jl}$ be the marginalizations of $w$ given by $w_j(x_j) = \int w(\mathbf{x})\,d\mathbf{x}_{-j}$ and $w_{jl}(x_j, x_l) = \int w(\mathbf{x})\,d\mathbf{x}_{-(j,l)}$. Define

$$\mathcal{H}(w) = \{\eta \in L_2(w) : \eta(\mathbf{x}) = \eta_1(x_1) + \cdots + \eta_d(x_d) \text{ for some functions}$$
$$\eta_1 \in L_2(w_1), \ldots, \eta_d \in L_2(w_d)\},$$

$$\mathcal{H}^0(w) = \left\{\eta \in \mathcal{H}(w) : \int \eta(\mathbf{x})w(\mathbf{x})\,d\mathbf{x} = 0\right\},$$

$$\mathcal{H}_j(w) = \{\eta \in \mathcal{H}(w) : \eta(\mathbf{x}) = \eta_j(x_j) \text{ for a function } \eta_j \in L_2(w_j)\},$$

$$\mathcal{H}_j^0(w) = \{\eta \in \mathcal{H}^0(w) : \eta(\mathbf{x}) = \eta_j(x_j) \text{ for a function } \eta_j \in L_2(w_j)\},$$

$$\mathcal{G}(w) = \{\boldsymbol{\eta} = (\eta_0, \eta_1, \ldots, \eta_d) : \eta_0 \in \mathbb{R} \text{ and } \eta_j \in \mathcal{H}_j(w) \text{ for } j = 1, \ldots, d\},$$

$$\mathcal{G}^0(w) = \{\boldsymbol{\eta} = (\eta_0, \eta_1, \ldots, \eta_d) : \eta_0 \in \mathbb{R} \text{ and } \eta_j \in \mathcal{H}_j^0(w) \text{ for } j = 1, \ldots, d\}.$$

The (semi)norm for functions $\eta \in \mathcal{H}(w)$ is defined by $\|\eta\|_w^2 = \int \eta^2(\mathbf{x}) w(\mathbf{x}) \, d\mathbf{x}$. For tuples of functions $\boldsymbol{\eta} \in \mathcal{G}(w)$ [or $\mathcal{G}^0(w)$], we define a Hilbert (semi)norm by $\|\boldsymbol{\eta}\|_w^2 = \int [\eta_0^2 + \sum_{j=1}^d \eta_j^2(x_j)] w(\mathbf{x}) \, d\mathbf{x}$. Within this framework, one can write

$$(7) \qquad \widehat{\mathbf{F}}\boldsymbol{\eta} = \widehat{\mathbf{F}}\boldsymbol{\eta}^0 + \widehat{\mathbf{F}}'(\boldsymbol{\eta}^0)(\boldsymbol{\eta} - \boldsymbol{\eta}^0) + o(\|\boldsymbol{\eta} - \boldsymbol{\eta}^0\|_{\widehat{w}^{\eta^0}}),$$

where $\widehat{\mathbf{F}}'(\boldsymbol{\eta}^0)(\cdot)$ is the Fréchet derivative of $\widehat{\mathbf{F}}$ at $\boldsymbol{\eta}^0$ in $L_2(\widehat{w}^{\eta^0})$ which is a linear transformation from $\mathcal{G}(\widehat{w}^{\eta^0})$ to $\mathcal{G}(\widehat{w}^{\eta^0})$. Its explicit form is given at (35) in Section 6.

The *outer loop* is originated from the linear approximation at (7). We adopt a Newton–Raphson iterative method for the outer loop. For simplicity, we write

$$\widehat{w}^{(k-1)} = \widehat{w}^{\widehat{\eta}^{(k-1)}},$$

$$\widehat{w}_j^{(k-1)}(x_j) = \int \widehat{w}^{(k-1)}(\mathbf{x}) \, d\mathbf{x}_{-j},$$

$$\widehat{w}_{jl}^{(k-1)}(x_j, x_l) = \int \widehat{w}^{(k-1)}(\mathbf{x}) \, d\mathbf{x}_{-(j,l)}.$$

Suppose that at the end of the $(k-1)$th outer iteration, or at the start of the $k$th outer iteration, we are given $\widehat{\boldsymbol{\eta}}^{(k-1)} = (\widehat{\eta}_0^{(k-1)}, \widehat{\eta}_1^{(k-1)}, \ldots, \widehat{\eta}_d^{(k-1)}) \in \mathcal{G}^0(\widehat{w}^{(k-1)})$. The updating equation for computing the $k$th outer iteration estimate is given by

$$(8) \qquad \mathbf{0} = \widehat{\mathbf{F}}\widehat{\boldsymbol{\eta}}^{(k-1)} + \widehat{\mathbf{F}}'(\widehat{\boldsymbol{\eta}}^{(k-1)})(\boldsymbol{\eta} - \widehat{\boldsymbol{\eta}}^{(k-1)}),$$

where $\widehat{\mathbf{F}}'(\widehat{\boldsymbol{\eta}}^{(k-1)})(\cdot)$ is the Fréchet derivative of $\widehat{\mathbf{F}}$ at $\widehat{\boldsymbol{\eta}}^{(k-1)}$, in $\mathcal{G}^0(\widehat{w}^{(k-1)})$. Define $\xi_j = \eta_j - \widehat{\eta}_j^{(k-1)}$, for $0 \le j \le d$, the changes in the $k$th outer update. The updating equation (8) can be written explicitly as the following system of equations:

$$\xi_0 = \left( \int \widehat{w}^{(k-1)}(\mathbf{x}) \, d\mathbf{x}, \right)^{-1}$$

$$(9) \qquad \times \int \left[ \frac{\widetilde{m}(\mathbf{x}) - g^{-1}(\widehat{\eta}^{(k-1)}(\mathbf{x}))}{V(g^{-1}(\widehat{\eta}^{(k-1)}(\mathbf{x}))) g'(g^{-1}(\widehat{\eta}^{(k-1)}(\mathbf{x})))} \right] \widehat{p}(\mathbf{x}) \, d\mathbf{x},$$

$$\xi_j(x_j) = \widetilde{\xi}_j^{(k)}(x_j) - \sum_{l=1, \neq j}^d \int \xi_l(x_l) \frac{\widehat{w}_{jl}^{(k-1)}(x_j, x_l)}{\widehat{w}_j^{(k-1)}(x_j)} \, dx_l - \xi_0, \qquad j = 1, \ldots, d,$$

where

$$\widetilde{\xi}_j^{(k)}(x_j) = \frac{\int \widetilde{\xi}^{(k)}(\mathbf{x}) \widehat{w}^{(k-1)}(\mathbf{x}) \, d\mathbf{x}_{-j}}{\int \widehat{w}^{(k-1)}(\mathbf{x}) \, d\mathbf{x}_{-j}}, \qquad j = 1, \ldots, d,$$

$$\widetilde{\xi}^{(k)}(\mathbf{x}) = \left[ \frac{\widetilde{m}(\mathbf{x}) - g^{-1}(\widehat{\eta}^{(k-1)}(\mathbf{x}))}{V(g^{-1}(\widehat{\eta}^{(k-1)}(\mathbf{x}))) g'(g^{-1}(\widehat{\eta}^{(k-1)}(\mathbf{x})))} \right] \frac{\widehat{p}(\mathbf{x})}{\widehat{w}^{(k-1)}(\mathbf{x})}.$$

The *inner loop* to get the $k$th outer iteration estimate is to find the solution $\xi_j, j = 0, \ldots, d$, of the system of equations (9). This is equivalent to finding the minimizer, in the space $\mathcal{H}(\widehat{w}^{(k-1)})$, of

$$\|\widetilde{\xi}^{(k)} - \xi\|^2_{\widehat{w}^{(k-1)}} = \int [\widetilde{\xi}^{(k)}(\mathbf{x}) - \xi_0 - \xi_1(x_1) - \cdots - \xi_d(x_d)]^2 \widehat{w}^{(k-1)}(\mathbf{x}) \, d\mathbf{x}$$

with the normalizing constraints $\int \xi_j(x_j) \widehat{w}^{(k-1)}(\mathbf{x}) \, d\mathbf{x} = 0, j = 1, \ldots, d$. The problem is exactly the same as the smooth backfitting of Mammen, Linton and Nielsen [15] except that the $L_2(\widehat{p})$ norm there is replaced by the $L_2(\widehat{w}^{(k-1)})$ norm. Thus, one can see that the smooth backfitting procedure based on (9) converges. Call the limit $\widehat{\xi}^{(k)}$. Note that $\widehat{\xi}^{(k)}(\mathbf{x})$ is uniquely decomposed into $\widehat{\xi}^{(k)}(\mathbf{x}) = \widehat{\xi}_0^{(k)} + \widehat{\xi}_1^{(k)}(x_1) + \cdots + \widehat{\xi}_d^{(k)}(x_d)$, where $\widehat{\xi}_0^{(k)} \in \mathbb{R}$ and $\widehat{\xi}_j^{(k)} \in \mathcal{H}_j^0(\widehat{w}^{(k-1)})$.

The components of the $k$th updated outer estimate are defined by

$$\widehat{\eta}_0^{(k)} = \widehat{\eta}_0^{(k-1)} + \widehat{\xi}_0^{(k)} + \sum_{j=1}^d c_j^{(k)},$$

(10)

$$\widehat{\eta}_j^{(k)}(x_j) = \widehat{\eta}_j^{(k-1)}(x_j) + \widehat{\xi}_j^{(k)}(x_j) - c_j^{(k)}, \qquad j = 1, \ldots, d,$$

where $c_j^{(k)} = [\int \widehat{w}_j^{(k)}(x_j) \, dx_j]^{-1} \int [\widehat{\eta}_j^{(k-1)}(x_j) + \widehat{\xi}_j^{(k)}(x_j)] \widehat{w}_j^{(k)} \, dx_j, j = 1, \ldots, d$. The tuple of these updated functions $\widehat{\boldsymbol{\eta}}^{(k)} = (\widehat{\eta}_0^{(k)}, \widehat{\eta}_1^{(k)}, \ldots, \widehat{\eta}_d^{(k)})$ equals the solution of the equation (8) in the space $\mathcal{G}^0(\widehat{w}^{(k)})$.

Returning to the inner loop, we note that the updating equation for the $j$th step of the $r$th iteration cycle is given by

$$\widehat{\xi}_j^{(k),[r]}(x_j) = \widetilde{\xi}_j^{(k)}(x_j) - \sum_{l<j} \int \widehat{\xi}_l^{(k),[r]}(x_l) \frac{\widehat{w}_{jl}^{(k-1)}(x_j, x_l)}{\widehat{w}_j^{(k-1)}(x_j)} \, dx_l$$

(11)

$$- \sum_{l>j} \int \widehat{\xi}_l^{(k),[r-1]}(x_l) \frac{\widehat{w}_{jl}^{(k-1)}(x_j, x_l)}{\widehat{w}_j^{(k-1)}(x_j)} \, dx_l - \widehat{\xi}_0^{(k)},$$

with $\widehat{\xi}_0^{(k)}$ defined by the first equation at (9). For an initial estimate in the inner iteration, one may take the centered version of $\widetilde{\xi}_j^{(k)}$: $\widehat{\xi}_j^{(k),[0]}(x_j) = \widetilde{\xi}_j^{(k)}(x_j) - \int \widetilde{\xi}_j^{(k)}(x_j) \widehat{w}_j^{(k-1)}(x_j) \, dx_j$. For an initial estimate $\widehat{\boldsymbol{\eta}}^{(0)}$ in the outer iteration, one may use some parametric model fits or use the marginal integration estimates.

**3. Estimation with local linear smoothing.** In this section, we propose maximum smoothed quasilikelihood estimation based on local linear fit. We briefly go over the projection interpretation of the local linear smooth backfitting in the ordinary additive models, which is the basic building block for the inner loop of our iterative algorithm. Here and in Section 4.2, we use the notation $\eta_0$, instead of $\eta$ in Section 2, to denote an additive function, and $\eta_j$, $1 \le j \le d$, to express its partial derivative with respect to $x_j$. The function $\eta_j$ does not mean the $j$th component of an additive function. For the latter, we write $\eta_{0j}$ instead.

### 3.1. *Projection property of local linear smoothers.*

To understand the full dimensional local linear fitting as a projection of the response vector $\mathbf{Y} = (Y^1, \ldots, Y^n)^T$ onto a relevant space, let the definitions of $\mathcal{F}$ and $\mathcal{F}_{\text{full}}$ in Section 2 be modified to $\mathcal{F} = \{(\mathbf{f}^1, \ldots, \mathbf{f}^n) : \mathbf{f}^i \in \mathcal{F}_0\}, \mathcal{F}_{\text{full}} = \{(\mathbf{f}, \ldots, \mathbf{f}) : \mathbf{f} \in \mathcal{F}\}$. Note that $\mathcal{F} = \mathcal{F}_0 \times \cdots \times \mathcal{F}_0$ and $\mathcal{F}_{\text{full}}$ is one-to-one correspondent to $\mathcal{F}_0$. The response vector can be embedded into $\mathcal{F}$ via $\mathbf{Y} \to (\mathbf{Y}^1, \ldots, \mathbf{Y}^n)$ where $\mathbf{Y}^i = (Y^i, 0, \ldots, 0)^T \in \mathcal{F}_0$.

Let $\mathbf{X}^i(\mathbf{x}) = (1, (X_1^i - x_1)/h_1, \ldots, (X_d^i - x_d)/h_d)^T$ and $K^i(\mathbf{x}) = n^{-1} K_\mathbf{h}(\mathbf{x}, \mathbf{X}^i)$. For a given $\mathbf{x}$, let $\widetilde{\beta}_0(\mathbf{x})$ be the full dimensional local linear estimator of $m(\mathbf{x})$, and $\widetilde{\beta}_j(\mathbf{x})$, for $1 \le j \le d$, be the full dimensional local linear estimator of $h_j \partial m(\mathbf{x})/\partial x_j$, respectively. Then, $\widetilde{\boldsymbol{\beta}}(\mathbf{x}) \equiv (\widetilde{\beta}_0(\mathbf{x}), \widetilde{\beta}_1(\mathbf{x}), \ldots, \widetilde{\beta}_d(\mathbf{x}))^T$ is given as the minimizer of the following quadratic form with respect to $\boldsymbol{\beta}(\mathbf{x}) = (\beta_0(\mathbf{x}), \beta_1(\mathbf{x}), \ldots, \beta_d(\mathbf{x}))^T$:

$$\sum_{i=1}^n [\mathbf{Y}^i - \boldsymbol{\beta}(\mathbf{x})]^T \mathbf{X}^i(\mathbf{x}) K^i(\mathbf{x}) \mathbf{X}^i(\mathbf{x})^T [\mathbf{Y}^i - \boldsymbol{\beta}(\mathbf{x})].$$

With the modified norm $\| \cdot \|_*$ defined by

$$\|(\mathbf{f}^1, \ldots, \mathbf{f}^n)\|_* = \left[ \int \sum_{i=1}^n \mathbf{f}^i(\mathbf{x})^T \mathbf{X}^i(\mathbf{x}) K^i(\mathbf{x}) \mathbf{X}^i(\mathbf{x})^T \mathbf{f}^i(\mathbf{x}) \, d\mathbf{x} \right]^{1/2},$$

the full dimensional estimator $\widetilde{\boldsymbol{\beta}}(\mathbf{x})$ can be regarded as a projection of $(\mathbf{Y}^1, \ldots, \mathbf{Y}^n)$ onto $\mathcal{F}_{\text{full}}$. It is also noted that for $(\mathbf{f}, \ldots, \mathbf{f}) \in \mathcal{F}_{\text{full}}$, the norm $\|(\mathbf{f}, \ldots, \mathbf{f})\|_*$ is simplified to $\|\mathbf{f}\|_{\widehat{V}} \equiv [\int \mathbf{f}(\mathbf{x})^T \widehat{\mathbf{V}}(\mathbf{x}) \mathbf{f}(\mathbf{x}) \, d\mathbf{x}]^{1/2}$ where $\widehat{\mathbf{V}}(\mathbf{x}) = \mathbf{X}(\mathbf{x})^T \mathbf{K}(\mathbf{x}) \mathbf{X}(\mathbf{x})$, and that $\| \cdot \|_{\widehat{V}}$ is an $L_2$-type norm for $\mathcal{F}_0$. For $(\boldsymbol{\beta}, \ldots, \boldsymbol{\beta}) \in \mathcal{F}_{\text{full}}$ with $\boldsymbol{\beta} \in \mathcal{F}_0$, the following Pythagorean identity holds:

$$\begin{aligned}
(12) \quad & \|(\mathbf{Y}^1, \ldots, \mathbf{Y}^n) - (\boldsymbol{\beta}, \ldots, \boldsymbol{\beta})\|_*^2 \\
& = \|(\mathbf{Y}^1, \ldots, \mathbf{Y}^n) - (\widetilde{\boldsymbol{\beta}}, \ldots, \widetilde{\boldsymbol{\beta}})\|_*^2 + \|\widetilde{\boldsymbol{\beta}} - \boldsymbol{\beta}\|_{\widehat{V}}^2.
\end{aligned}$$

The identity (12) suggests a clue to construct an estimator for a structured model. If one assumes a model class which is a subspace of $\mathcal{F}_0$, then one can get an $M$-type estimator by minimizing the second term on the right-hand side of (12) over the assumed model class. For a matrix-valued function $\mathbf{V}$ for which $\mathbf{V}(\mathbf{x})$ is positive definite for all $\mathbf{x}$, let $\mathcal{F}_0(\mathbf{V})$ denote the space $\mathcal{F}_0$ equipped with the norm $\|\mathbf{f}\|_{\mathbf{V}} = [\int \mathbf{f}(\mathbf{x})^T \mathbf{V}(\mathbf{x}) \mathbf{f}(\mathbf{x}) \, d\mathbf{x}]^{1/2}$. The definition of the space $\mathcal{H}$ in Section 2 is modified to

$$\mathcal{H}(\mathbf{V}) = \{\mathbf{f} \in \mathcal{F}_0(\mathbf{V}) : f_0(\mathbf{x}) = f_{01}(x_1) + \cdots + f_{0d}(x_d) \text{ for some functions}$$

$$f_{0j} : \mathbb{R} \to \mathbb{R}, \text{ and } f_j(\mathbf{x}) = g_j(x_j) \text{ for some function}$$

$$g_j : \mathbb{R} \to \mathbb{R}, j = 1, \ldots, d\}.$$

Then, the local linear smooth backfitting estimator in the ordinary additive models, proposed by Mammen, Linton and Nielsen [15], can be given as the projection of the full dimensional local linear estimator $\widetilde{\boldsymbol{\beta}}$ onto $\mathcal{H}(\widehat{\mathbf{V}})$.

For $j = 1, \ldots, d$, define

$$\mathcal{H}_j(\mathbf{V}) = \{\mathbf{f} \in \mathcal{H}(\mathbf{V}) : f_0(\mathbf{x}) = f_{0j}(x_j), f_k \equiv 0 \text{ for } k \neq j\}.$$

The space $\mathcal{H}(\mathbf{V})$ equals $\mathcal{H}_1(\mathbf{V}) + \cdots + \mathcal{H}_d(\mathbf{V})$. Let $\Pi_{j,\mathbf{V}}$ denote the projection operator onto $\mathcal{H}_j(\mathbf{V})$. To express the projections explicitly, let $\mathbf{M}_{j,\mathbf{V}}(x_j)$ be a $2 \times 2$ matrix and $\mathbf{A}_j$ be a $2 \times (d+1)$ matrix such that

$$(13) \qquad \mathbf{M}_{j,\mathbf{V}}(x_j) = \begin{bmatrix} V_{00,j}(x_j) & V_{0j,j}(x_j) \\ V_{0j,j}(x_j) & V_{jj,j}(x_j) \end{bmatrix} \quad \text{and} \quad \mathbf{A}_j = \begin{bmatrix} \mathbf{1}_0^T \\ \mathbf{1}_j^T \end{bmatrix},$$

where $V_{pq,j}(x_j)$ are $(p,q)$th elements of the matrix $\mathbf{V}_j(x_j) \equiv \int \mathbf{V}(\mathbf{x}) \, d\mathbf{x}_{-j}$, and $\mathbf{1}_k$ is a $(d+1)$-dimensional unit vector with 1 appearing at the $(k+1)$th position. Then, it can be shown that for $\mathbf{f} \in \mathcal{H}(\mathbf{V})$,

$$(\Pi_{j,\mathbf{V}}\mathbf{f})(x_j) = (g_{0j}(x_j), 0, \ldots, 0, g_j(x_j), 0, \ldots, 0)^T$$

where

$$(g_{0j}(x_j), g_j(x_j))^T = \mathbf{M}_{j,\mathbf{V}}(x_j)^{-1} \int \mathbf{A}_j \mathbf{V}(\mathbf{x}) \mathbf{f}(\mathbf{x}) \, d\mathbf{x}_{-j}.$$

Since $\mathbf{1}_0 \in \mathcal{H}_j(\mathbf{V})$ for all $j = 1, \ldots, d$, the decomposition of $\mathbf{f} \in \mathcal{H}(\mathbf{V})$ into $\mathbf{f}(\mathbf{x}) = \mathbf{f}_1(\mathbf{x}) + \cdots + \mathbf{f}_d(\mathbf{x})$ with $\mathbf{f}_j \in \mathcal{H}_j(\mathbf{V})$ is not unique. For a unique identification, let

$$\mathcal{H}_j^0(\mathbf{V}) = \{\mathbf{f} \in \mathcal{H}(\mathbf{V}) : f_0(\mathbf{x}) = f_{0j}(x_j), f_k \equiv 0 \text{ for } k \neq j, \langle \mathbf{f}, \mathbf{1}_0 \rangle_{\mathbf{V}} = 0\},$$

where $\langle \mathbf{f}, \mathbf{g} \rangle_{\mathbf{V}} = \int \mathbf{f}^T(\mathbf{x}) \mathbf{V}(\mathbf{x}) \mathbf{g}(\mathbf{x}) \, d\mathbf{x}$. The norming constraint $\langle \mathbf{f}, \mathbf{1}_0 \rangle_{\mathbf{V}} = 0$ implies that $\mathbf{f}$ is orthogonal to constant functions, which is equivalent to the centering constraint in the local constant case. The local linear smooth backfitting estimator $\widehat{\boldsymbol{\beta}}$ in the ordinary additive models can be written as $\widehat{\boldsymbol{\beta}}(\mathbf{x}) = \widehat{\boldsymbol{\beta}}_0 + \widehat{\boldsymbol{\beta}}_1(x_1) + \cdots + \widehat{\boldsymbol{\beta}}_d(x_d)$ where $\widehat{\boldsymbol{\beta}}_0 = \overline{Y}\mathbf{1}_0$ and $\widehat{\boldsymbol{\beta}}_j$ $(j = 1, \ldots, d)$ satisfy the following system of linear integral equations:

$$(14) \qquad \begin{aligned} \widehat{\boldsymbol{\beta}}_j &= \widetilde{\boldsymbol{\beta}}_j - \sum_{l=1, \neq j}^{d} \Pi_{j,\widehat{\mathbf{V}}}(\widehat{\boldsymbol{\beta}}_l) - \widehat{\boldsymbol{\beta}}_0, \qquad j = 1, \ldots, d, \\ \langle \widehat{\boldsymbol{\beta}}_j, \mathbf{1}_0 \rangle_{\widehat{\mathbf{V}}} &= 0, \qquad\qquad\qquad\qquad\qquad j = 1, \ldots, d. \end{aligned}$$

Here, $\widetilde{\boldsymbol{\beta}}_j(x_j) = (\widetilde{\beta}_{0j}(x_j), 0, \dots, 0, \widetilde{\beta}_j(x_j), 0, \dots, 0)^T$ and $(\widetilde{\beta}_{0j}(x_j), \widetilde{\beta}_j(x_j))^T$ denotes the vector of the marginal local linear estimators of $E(Y^1 | X_j^1 = x_j)$ and its derivative, obtained by regressing $Y^i$ on $X_j^i$ only. The local linear smooth backfitting estimator of $m_j(x_j)$, in $m(\mathbf{x}) = m_0 + m_1(x_1) + \cdots + m_d(x_d)$ with $Em_j(X_j^1) = 0$ for $1 \le j \le d$, equals $\widehat{\beta}_{0j}(x_j)$, and that of its derivative $\partial m_j(x_j)/\partial x_j$ equals $\widehat{\beta}_j(x_j)/h_j$.

3.2. *The smoothed quasilikelihood and backfitting algorithms.* In this subsection, we let $\eta_0^*(\mathbf{x}) = \eta_{00}^* + \eta_{01}^*(x_1) + \cdots + \eta_{0d}^*(x_d)$ denote the true additive function, where each component $\eta_{0j}^*$ is defined by (5). Also, let $\eta_j^*(x_j) = h_j \eta_{0j}^{*\prime}(x_j)$ for $1 \le j \le d$. The function $\eta_j^*$ should not be confused with $\eta_{0j}^*$, the $j$th component function of $\eta_0^*$. They are the targets of the maximum smoothed quasilikelihood estimators $\widehat{\eta}_0, \widehat{\eta}_1, \dots, \widehat{\eta}_d$ that we describe below.

For $\boldsymbol{\eta} = (\eta_0, \eta_1, \dots, \eta_d)^T \in \mathcal{F}_0$, define

$$\boldsymbol{\eta}(\mathbf{u}, \mathbf{x}) = \eta_0(\mathbf{x}) + \left(\frac{u_1 - x_1}{h_1}\right)\eta_1(\mathbf{x}) + \cdots + \left(\frac{u_d - x_d}{h_d}\right)\eta_d(\mathbf{x}).$$

We include $\eta_j(\mathbf{x})$ for $1 \le j \le d$ in $\boldsymbol{\eta}(\mathbf{x})$ to put the problems of estimating $\eta_0$ and its derivatives into the same framework of projection operation. With a general link $g$, we define a *smoothed quasilikelihood* for local linear fit by

$$SQ(\boldsymbol{\eta}) = \int n^{-1} \sum_{i=1}^n Q(g^{-1}(\boldsymbol{\eta}(\mathbf{X}^i, \mathbf{x})), Y^i) K_{\mathbf{h}}(\mathbf{x}, \mathbf{X}^i)\, d\mathbf{x}.$$

We call $\boldsymbol{\eta}$ additive if $\boldsymbol{\eta} \in \mathcal{H}$, that is, $\eta_0(\mathbf{x}) = \eta_{00} + \eta_{01}(x_1) + \cdots + \eta_{0d}(x_d)$ and $\eta_j(\mathbf{x}) = \eta_j(x_j)$, $j = 1, \dots, d$. We define $\widehat{\boldsymbol{\eta}}$ to be the maximizer of the smoothed quasilikelihood $SQ(\boldsymbol{\eta})$ over all additive functions $\boldsymbol{\eta}$. Each additive function $\boldsymbol{\eta}$ can be written as

$$\boldsymbol{\eta} = \boldsymbol{\eta}_0 + \boldsymbol{\eta}_1 + \cdots + \boldsymbol{\eta}_d$$

where $\boldsymbol{\eta}_0 = \eta_{00}\mathbf{1}_0$ and $\boldsymbol{\eta}_j(\mathbf{x}) = (\eta_{0j}(x_j), 0, \dots, 0, \eta_j(x_j), 0, \dots, 0)^T$. We consider the following space:

$$\mathcal{G}^0(\mathbf{V}) = \{(\boldsymbol{\eta}_0, \boldsymbol{\eta}_1, \dots, \boldsymbol{\eta}_d) : \boldsymbol{\eta}_0 = \eta_{00}\mathbf{1}_0 \text{ for } \eta_{00} \in \mathbb{R}, \boldsymbol{\eta}_j \in \mathcal{H}_j^0(\mathbf{V}), j = 1, \dots, d\}.$$

The space $\mathcal{G}^0$ is endowed with a Hilbert (semi) norm defined by

$$\|(\boldsymbol{\eta}_0, \boldsymbol{\eta}_1, \dots, \boldsymbol{\eta}_d)\|_{\mathbf{V}}^2$$
$$= |\boldsymbol{\eta}_0|^2 \int V_{00}(\mathbf{x})\, d\mathbf{x} + \sum_{j=1}^d \int \boldsymbol{\eta}_j(x_j)^T \left(\int \mathbf{V}(\mathbf{x})\, d\mathbf{x}_{-j}\right)\boldsymbol{\eta}_j(x_j)\, dx_j.$$

With a slight abuse of notation we continue to use $\|\cdot\|_{\mathbf{V}}$ for the norm of $\mathcal{G}^0$ as we use it for the norm of $\mathcal{H}$.

Let $\boldsymbol{\eta}_L$ denote the element of $\mathcal{G}^0$ that corresponds to an additive function $\boldsymbol{\eta} \in \mathcal{H}$. With this convention, define

$$\widehat{F}_{00}\boldsymbol{\eta}_L = \int n^{-1} \sum_{i=1}^{n} q_1(\boldsymbol{\eta}(\mathbf{X}^i, \mathbf{x}), Y^i) K_{\mathbf{h}}(\mathbf{x}, \mathbf{X}^i) \, d\mathbf{x},$$

$$\widehat{F}_{0j}\boldsymbol{\eta}_L = \int n^{-1} \sum_{i=1}^{n} q_1(\boldsymbol{\eta}(\mathbf{X}^i, \mathbf{x}), Y^i) K_{\mathbf{h}}(\mathbf{x}, \mathbf{X}^i) \, d\mathbf{x}_{-j}, \qquad j = 1, \ldots, d,$$

$$\widehat{F}_{j}\boldsymbol{\eta}_L = \int n^{-1} \sum_{i=1}^{n} q_1(\boldsymbol{\eta}(\mathbf{X}^i, \mathbf{x}), Y^i) \left( \frac{X_j^i - x_j}{h_j} \right) K_{\mathbf{h}}(\mathbf{x}, \mathbf{X}^i) \, d\mathbf{x}_{-j},$$
$$j = 1, \ldots, d,$$

$$(\widehat{\mathbf{F}}\boldsymbol{\eta}_L)(\mathbf{x}) = (\widehat{F}_{00}\boldsymbol{\eta}_L, (\widehat{F}_{01}\boldsymbol{\eta}_L)(x_1), \ldots, (\widehat{F}_{0d}\boldsymbol{\eta}_L)(x_d),$$
$$(\widehat{F}_{1}\boldsymbol{\eta}_L)(x_1), \ldots, (\widehat{F}_{d}\boldsymbol{\eta}_L)(x_d))^T.$$

Then, $\widehat{\boldsymbol{\eta}}_L$ that corresponds to $\widehat{\boldsymbol{\eta}}$ may be obtained by solving $\widehat{\mathbf{F}}\boldsymbol{\eta}_L = \mathbf{0}$ for $\boldsymbol{\eta}_L \in \mathcal{G}^0$. As in Section 2, we approximate $\widehat{\mathbf{F}}\boldsymbol{\eta}_L$ for $\boldsymbol{\eta}_L$ in a neighborhood of $\boldsymbol{\eta}_L^0$. To do this we need to consider a proper metric for $\mathcal{G}^0$. Define

$$\widehat{w}^i(\mathbf{x}, \boldsymbol{\eta}) = -q_2(\boldsymbol{\eta}(\mathbf{X}^i, \mathbf{x}), Y^i) K_{\mathbf{h}}(\mathbf{x}, \mathbf{X}^i),$$
$$\widehat{\mathbf{V}}(\mathbf{x}, \boldsymbol{\eta}) = \mathbf{X}(\mathbf{x})^T (n^{-1} \mathrm{diag}[\widehat{w}^1(\mathbf{x}, \boldsymbol{\eta}), \ldots, \widehat{w}^n(\mathbf{x}, \boldsymbol{\eta})]) \mathbf{X}(\mathbf{x}).$$

Then, writing $\widehat{\mathbf{V}}^{(0)} = \widehat{\mathbf{V}}(\mathbf{x}, \boldsymbol{\eta}^0)$ we have

$$(15) \qquad \widehat{\mathbf{F}}\boldsymbol{\eta}_L = \widehat{\mathbf{F}}\boldsymbol{\eta}_L^0 + \widehat{\mathbf{F}}'(\boldsymbol{\eta}_L^0)(\boldsymbol{\eta}_L - \boldsymbol{\eta}_L^0) + o(\|\boldsymbol{\eta}_L - \boldsymbol{\eta}_L^0\|_{\widehat{\mathbf{V}}^{(0)}}),$$

where $\widehat{\mathbf{F}}'(\boldsymbol{\eta}_L^0)(\cdot)$ is the Fréchet derivative of $\widehat{\mathbf{F}}$ at $\boldsymbol{\eta}_L^0$ in $\mathcal{G}^0(\widehat{\mathbf{V}}^{(0)})$.

As in Section 2, the *outer loop* for solving $\widehat{\mathbf{F}}\boldsymbol{\eta}_L = \mathbf{0}$ can be based on the linear approximation (15). The updating equation for computing the $k$th outer iteration estimate $\widehat{\boldsymbol{\eta}}_L^{(k)}$ is given by

$$(16) \qquad \mathbf{0} = \widehat{\mathbf{F}}\widehat{\boldsymbol{\eta}}_L^{(k-1)} + \widehat{\mathbf{F}}'(\widehat{\boldsymbol{\eta}}_L^{(k-1)})(\boldsymbol{\eta}_L - \widehat{\boldsymbol{\eta}}_L^{(k-1)}),$$

where $\widehat{\mathbf{F}}'(\widehat{\boldsymbol{\eta}}_L^{(k-1)})(\cdot)$ is the Fréchet derivative of $\widehat{\mathbf{F}}$ at $\widehat{\boldsymbol{\eta}}_L^{(k-1)}$ in $\mathcal{G}^0(\widehat{\mathbf{V}}^{(k-1)})$ and $\widehat{\mathbf{V}}^{(k-1)} = \widehat{\mathbf{V}}(\mathbf{x}, \widehat{\boldsymbol{\eta}}^{(k-1)})$. Let $\xi_{00} = \eta_{00} - \widehat{\eta}_{00}^{(k-1)}$, $\xi_{0j} = \eta_{0j} - \widehat{\eta}_{0j}^{(k-1)}$ and $\xi_j = \eta_j - \widehat{\eta}_j^{(k-1)}$. To get an explicit form of the updating equation (16), define $\widehat{\mathbf{M}}_j^{(k-1)}(x_j) \equiv \mathbf{M}_{j,\widehat{\mathbf{V}}^{(k-1)}}(x_j)$ in the same way as $\mathbf{M}_{j,\mathbf{V}}(x_j)$ at (13) with $\mathbf{V}$ replaced by $\widehat{\mathbf{V}}^{(k-1)}$. Also, define

$$\widehat{\mathbf{M}}_{j,l}^{(k-1)}(x_j, x_l) = \begin{bmatrix} \widehat{V}_{00,jl}^{(k-1)}(x_j, x_l) & \widehat{V}_{0l,jl}^{(k-1)}(x_j, x_l) \\ \widehat{V}_{0j,jl}^{(k-1)}(x_j, x_l) & \widehat{V}_{jl,jl}^{(k-1)}(x_j, x_l) \end{bmatrix},$$

where $\widehat{V}_{pq,jl}^{(k-1)}(x_j, x_l)$ are $(p,q)$th elements of the matrix $\widehat{\mathbf{V}}_{jl}^{(k-1)}(x_j, x_l) \equiv \int \widehat{\mathbf{V}}^{(k-1)}(\mathbf{x})\, d\mathbf{x}_{-(j,l)}$. Furthermore, for $j = 1, \ldots, n$ we let

$$\widetilde{\zeta}_{0j}^{(k)}(x_j) = -\int \frac{1}{n} \sum_{i=1}^{n} \left[ \frac{q_1(\widehat{\boldsymbol{\eta}}^{(k-1)}(\mathbf{X}^i, \mathbf{x}), Y^i)}{q_2(\widehat{\boldsymbol{\eta}}^{(k-1)}(\mathbf{X}^i, \mathbf{x}), Y^i)} \right] \widehat{w}^i(\mathbf{x}, \widehat{\boldsymbol{\eta}}^{(k-1)})\, d\mathbf{x}_{-j},$$

$$\widetilde{\zeta}_{j}^{(k)}(x_j) = -\int \frac{1}{n} \sum_{i=1}^{n} \left[ \frac{q_1(\widehat{\boldsymbol{\eta}}^{(k-1)}(\mathbf{X}^i, \mathbf{x}), Y^i)}{q_2(\widehat{\boldsymbol{\eta}}^{(k-1)}(\mathbf{X}^i, \mathbf{x}), Y^i)} \right] \left( \frac{X_j^i - x_j}{h_j} \right) \widehat{w}^i(\mathbf{x}, \widehat{\boldsymbol{\eta}}^{(k-1)})\, d\mathbf{x}_{-j}.$$

Then, it can be shown that the updating equation (16) is equivalent to

$$
\begin{aligned}
\widehat{\mathbf{M}}_j^{(k-1)}(x_j) \begin{bmatrix} \xi_{0j}(x_j) \\ \xi_j(x_j) \end{bmatrix} &= \begin{bmatrix} \widetilde{\zeta}_{0j}^{(k)}(x_j) \\ \widetilde{\zeta}_{j}^{(k)}(x_j) \end{bmatrix} - \xi_{00} \begin{bmatrix} \widehat{V}_{00,j}^{(k-1)}(x_j) \\ \widehat{V}_{0j,j}^{(k-1)}(x_j) \end{bmatrix} \\
&\quad - \sum_{l=1, \neq j}^{d} \int \widehat{\mathbf{M}}_{j,l}^{(k-1)}(x_j, x_l) \begin{bmatrix} \xi_{0l}(x_l) \\ \xi_l(x_l) \end{bmatrix} dx_l,
\end{aligned}
$$

(17)

$$
\begin{aligned}
\xi_{00} &= \left[ -\int \frac{1}{n} \sum_{i=1}^{n} \left[ \frac{q_1(\widehat{\boldsymbol{\eta}}^{(k-1)}(\mathbf{X}^i, \mathbf{x}), Y^i)}{q_2(\widehat{\boldsymbol{\eta}}^{(k-1)}(\mathbf{X}^i, \mathbf{x}), Y^i)} \right] \widehat{w}^i(\mathbf{x}, \widehat{\boldsymbol{\eta}}^{(k-1)})\, d\mathbf{x} \right] \\
&\quad \times \left[ \int n^{-1} \sum_{i=1}^{n} \widehat{w}^i(\mathbf{x}, \widehat{\boldsymbol{\eta}}^{(k-1)})\, d\mathbf{x} \right]^{-1},
\end{aligned}
$$

with the normalizing constraint

$$\int n^{-1} \sum_{i=1}^{n} \widehat{w}^i(\mathbf{x}, \widehat{\boldsymbol{\eta}}^{(k-1)}) \left[ \xi_{0j}(x_j) + \left( \frac{X_j^i - x_j}{h_j} \right) \xi_j(x_j) \right] d\mathbf{x} = 0, \qquad j = 1, \ldots, d.$$

Solving (17) constitutes our *inner loop* to find the $k$th outer iteration estimate. The system of equations (17) can be written in a different form using a projection operator as in (14). To do this, define

$$[\widetilde{\xi}_{0j}^{(k)}(x_j), \widetilde{\xi}_{j}^{(k)}(x_j)]^T = \widehat{\mathbf{M}}_j^{(k-1)}(x_j)^{-1} [\widetilde{\zeta}_{0j}^{(k)}(x_j), \widetilde{\zeta}_{j}^{(k)}(x_j)]^T$$

and $\widetilde{\boldsymbol{\xi}}_j^{(k)}(x_j) = (\widetilde{\xi}_{0j}^{(k)}(x_j), 0, \ldots, 0, \widetilde{\xi}_{j}^{(k)}(x_j), 0, \ldots, 0)^T \in \mathcal{H}_j$. Write

$$\boldsymbol{\xi}_0 = \xi_{00} \mathbf{1}_0,$$
$$\boldsymbol{\xi}_j(x_j) = (\xi_{0j}(x_j), 0, \ldots, 0, \xi_j(x_j), 0, \ldots, 0)^T, \qquad j = 1, \ldots, d.$$

Let $\widehat{\Pi}_j^{(k-1)} \equiv \Pi_{j, \widehat{\mathbf{V}}^{(k-1)}}$ be the projection operator onto $\mathcal{H}_j(\widehat{\mathbf{V}}^{(k-1)})$. Then, solving (17) is equivalent to solving

(18) $$\boldsymbol{\xi}_j = \widetilde{\boldsymbol{\xi}}_j^{(k)} - \sum_{l=1, \neq j}^{d} \widehat{\Pi}_j^{(k-1)}(\boldsymbol{\xi}_l) - \boldsymbol{\xi}_0, \qquad j = 1, \ldots, d,$$

subject to the normalizing constraint $\langle \boldsymbol{\xi}_j, \mathbf{1}_0 \rangle_{\widehat{\mathbf{V}}^{(k-1)}} = 0$, $j = 1, \ldots, d$.

The smooth backfitting algorithm based on (18) converges since it has the same projection interpretation as the local linear smooth backfitting in ordinary additive regression. Let $\widehat{\zeta}_{00}^{(k)}, \widehat{\zeta}_{0j}^{(k)}, \widehat{\zeta}_j^{(k)}$ denote the solution of the system of equations (18). Define for $j = 1, \ldots, d$

$$
\begin{aligned}
c_j^{(k)} = {}& \left[ \int n^{-1} \sum_{i=1}^n \widehat{w}^i(\mathbf{x}, \widehat{\boldsymbol{\eta}}^{(k)}) \, d\mathbf{x} \right]^{-1} \\
& \times \int n^{-1} \sum_{i=1}^n \widehat{w}^i(\mathbf{x}, \widehat{\boldsymbol{\eta}}^{(k)}) \Big[ \eta_{0j}^{(k-1)}(x_j) + \xi_{0j}^{(k)}(x_j) \\
& \qquad\qquad + \left( \frac{X_j^i - x_j}{h_j} \right) (\eta_j^{(k-1)}(x_j) + \xi_j^{(k)}(x_j)) \Big] \, d\mathbf{x}.
\end{aligned}
$$

Then, the $k$th outer iteration updates are given by

$$
\begin{aligned}
\widehat{\eta}_{00}^{(k)} &= \widehat{\eta}_{00}^{(k-1)} + \widehat{\xi}_{00}^{(k)} + \sum_{j=1}^d c_j^{(k)}, \\
\widehat{\eta}_{0j}^{(k)}(x_j) &= \widehat{\eta}_{0j}^{(k-1)}(x_j) + \widehat{\xi}_{0j}^{(k)}(x_j) - c_j^{(k)}, \qquad j = 1, \ldots, d, \\
\widehat{\eta}_j^{(k)}(x_j) &= \widehat{\eta}_j^{(k-1)}(x_j) + \widehat{\xi}_j^{(k)}(x_j), \qquad\qquad j = 1, \ldots, d.
\end{aligned}
$$

**4. Asymptotic and algorithmic properties.** First, we collect the assumptions for the theoretical results to be presented in this section.

ASSUMPTIONS.

A1. $p$ is bounded away from zero and infinity on its support, $[0, 1]^d$, and has continuous partial derivatives.

A2. $q_2(u, y) < 0$ for $u \in \mathbb{R}$ and $y$ in the range of the response, the link $g$ is strictly monotone and is three times continuously differentiable, $V$ is strictly positive and twice continuously differentiable, and $v(\mathbf{x}) \equiv \mathrm{var}(Y | \mathbf{X} = \mathbf{x})$ is continuous. $E|Y|^{r_0} < \infty$ for some $r_0 > 5/2$.

A3. The true component functions $\eta_j^*$'s in Section 2 and $\eta_{0j}^*$ in Section 3 are twice continuously differentiable.

A4. The base kernel function $K^0$ is a symmetric density function with compact support, $[-1, 1]$ say, and is Lipschitz continuous.

A5. $n^{1/5} h_j$ converge to constants $\delta_j > 0$ for $j = 1, \ldots, d$ as $n$ goes to infinity.

4.1. *Nadaraya–Watson smooth backfitting.* The first two theorems are for the limiting distributions of $\widehat{\eta}_j(x_j)$, $j = 1, \ldots, d$, defined by (6).

THEOREM 1 (Rates of convergence). *Suppose that the conditions* A1–A5 *hold. Then*

$$\|\hat{\eta} - \eta^*\|_p = O_p(n^{-2/5}),$$

$$\sup_{x_j \in [2h_j, 1-2h_j]} |\hat{\eta}_j(x_j) - \eta_j^*(x_j)| = O_p(n^{-2/5}\sqrt{\log n}), \qquad j = 1, \ldots, d.$$

For the statement of the next theorem, let $v(\mathbf{x}) = \mathrm{Var}(Y|\mathbf{X} = \mathbf{x})$, and write for simplicity $w^*(\mathbf{x}) = w^{\eta^*}(\mathbf{x})$. Define, for $\delta_j$ in the condition A5,

$$
\begin{aligned}
(19) \quad v_j(x_j) &= \frac{E[v(\mathbf{X})V(g^{-1}(\eta^*(\mathbf{X})))^{-2}g'(g^{-1}(\eta^*(\mathbf{X})))^{-2}|X_j = x_j]}{E^2[V(g^{-1}(\eta^*(\mathbf{X})))^{-1}g'(g^{-1}(\eta^*(\mathbf{X})))^{-2}|X_j = x_j]} \\
&\quad \times \delta_j^{-1}p_j(x_j)^{-1}\int [K^0(t)]^2\,dt,
\end{aligned}
$$

$$
\begin{aligned}
(20) \quad \beta(\mathbf{x}) &= -\sum_{j=1}^d \delta_j^2 \left[ p^{-1}(\mathbf{x})\frac{\partial}{\partial x_j}p(\mathbf{x})\frac{\partial}{\partial x_j}g^{-1}(\eta^*(\mathbf{x})) + \frac{1}{2}\frac{\partial^2}{\partial x_j^2}g^{-1}(\eta^*(\mathbf{x})) \right] \\
&\quad \times g'(g^{-1}(\eta^*(\mathbf{x})))\int t^2 K^0(t)\,dt.
\end{aligned}
$$

Let the constant $b_0$ and the functions $\beta_j(x_j)$ minimize $\int[\beta(\mathbf{x}) - b_0 - \sum_{j=1}^d \beta_j(x_j)]^2 w^*(\mathbf{x})\,d\mathbf{x}$, subject to $\int \beta_j(x_j)w_j^*(x_j)\,dx_j = 0$ for $j = 1, \ldots, d$.

THEOREM 2 (Asymptotic distributions). *Under the conditions of Theorem* 1, *for any* $x_1, \ldots, x_d \in (0, 1)$, $n^{2/5}[\hat{\eta}_1(x_1) - \eta_1^*(x_1), \ldots, \hat{\eta}_d(x_d) - \eta_d^*(x_d)]^T$ *converges in distribution to the* $d$-*variate normal distribution with mean vector* $[\beta_1(x_1), \ldots, \beta_d(x_d)]^T$ *and variance matrix* $\mathrm{diag}[v_j(x_j)]$.

Unlike the smooth backfitting estimator in the ordinary additive models, our estimator of the intercept $\eta_0^*$ has a nonnegligible asymptotic bias. In fact,

$$n^{2/5}(\hat{\eta}_0 - \eta_0^*) \xrightarrow{p} \beta_0,$$

where $\beta_0$ has a complicated form and is different from $b_0$ defined above. Writing $\mu_j = \int_{-1}^1 u^j K^0(u)\,du$ and $\kappa = \int_0^1 [\mu_1(-t)/\mu_0(-t)]\,dt$ where $\mu_j(c) = \int_c^1 u^j K^0(u)\,du$, it can be shown

$$
\begin{aligned}
(21) \quad \beta_0 &= E(q_2(\eta^*(\mathbf{X}^1), g^{-1}(\eta^*(\mathbf{X}^1))))^{-1} \\
&\quad \times \sum_{j=1}^d \delta_j^2 \left[ \tfrac{1}{2}\mu_2 \int \varphi_{jj}(\mathbf{x})\omega(\mathbf{x})\,d\mathbf{x} \right.
\end{aligned}
$$

$$+ \kappa \int \{\varphi_j(0+, \mathbf{x}_{-j}) \omega(0+, \mathbf{x}_{-j})$$
$$- \varphi_j(1-, \mathbf{x}_{-j}) \omega(1-, \mathbf{x}_{-j})\} \, d\mathbf{x}_{-j}\Big],$$

where $\varphi_j(a, \mathbf{x}_{-j}) = (\partial/\partial x_j)g^{-1}(\eta^*(\mathbf{x}))$, $\varphi_{jj}(\mathbf{x}) = (\partial^2/\partial x_j^2)g^{-1}(\eta^*(\mathbf{x}))$ and $\omega(\mathbf{x}) = w^*(\mathbf{x}) \times g'(g^{-1}(\eta^*(\mathbf{x})))$. Here, the argument $(a, \mathbf{x}_{-j})$ implies $\mathbf{x}$ with $x_j$ being replaced by $a$. From Theorem 2 and the convergence of $\widehat{\eta}_0$, we have, for $\mathbf{x}$ in the interior of the support of $p$,

$$n^{2/5}(\widehat{\eta}(\mathbf{x}) - \eta^*(\mathbf{x})) \quad \overset{d}{\Longrightarrow} \quad N\left(\beta_0 + \sum_{j=1}^{d} \beta_j(x_j), \sum_{j=1}^{d} v_j(x_j)\right).$$

Theorems 1 and 2 show that the proposed estimator has the desirable dimension reduction property. It achieves the same convergence rates as one-dimensional estimators. Furthermore, the asymptotic variance of $\widehat{\eta}_j(x_j)$ coincides with that of the one-dimensional local constant estimator obtained by fitting the model $E(Y|\mathbf{X} = \mathbf{x}) = g^{-1}(\eta_j(x_j) + \sum_{k=1, \neq j}^{d} \eta_k^*(x_k))$ with the other component functions $\eta_k^*$ ($k \neq j$) being known, see Fan, Heckman and Wand [5], for example. In this sense, our estimator $\widehat{\eta}_j(x_j)$ of the $j$th component function $\eta_j(x_j)$ enjoys the oracle variance.

REMARK 1. Theorems 1 and 2 hold regardless of whether or not $V$ correctly models the conditional variance of the response variable.

REMARK 2. Simultaneous confidence intervals for $\eta_j^*$ may be constructed using the joint limit distribution given in Theorem 2. This would involve estimation of $\beta_j$ and $v_j$ which is typically harder than the original problem of estimating $\eta^*$. Instead, one may use a bootstrap method.

REMARK 3. In the case where $Q(m, y) = -(y - m)^2/2$ and the link $g$ is the identity function, our results coincide with those of Mammen, Linton and Nielsen [15]. In this sense, our maximum smoothed quasilikelihood estimator can be regarded as an extension of the smooth backfitting to the context of generalized additive models.

The next two theorems are for the convergence of the proposed outer and inner iterative algorithms. Note that the uniform convergence of the inner iteration in Theorem 4 is required for the entire iteration to converge. Let $B_r(\widehat{\boldsymbol{\eta}})$ denote the ball centered at $\widehat{\boldsymbol{\eta}}$ with a radius $r$.

THEOREM 3 (Convergence of outer iteration). *Let $\widehat{\boldsymbol{\eta}}^{(k)}$ be the $k$th outer step estimator defined by* (8). *Under conditions* A1–A5, *there exist fixed* *r, C > 0 and 0 < γ < 1 which have the following property: if the initial estimator $\widehat{\boldsymbol{\eta}}^{(0)}$ belongs to $B_r(\widehat{\boldsymbol{\eta}})$ with probability tending to one, then*

$$\|\widehat{\boldsymbol{\eta}}^{(k)} - \widehat{\boldsymbol{\eta}}\|_p \leq C2^{-(k-1)}\gamma^{2^k-1}$$

*with probability tending to one.*

THEOREM 4 (Convergence of inner iteration). *Under conditions* A1–A5, *the inner iteration converges at a geometric rate. Moreover, if the initial estimator belongs to the ball introduced in Theorem* 3 *with probability tending to one, then the geometric convergence of the inner iteration is uniform for all steps in the outer iteration, with probability tending to one.*

REMARK 4. In practice, a parametric model fit can be used as an initial estimator. In our numerical experiments, the maximum likelihood estimator of the constant model, $\widehat{\eta}_0^{(0)} = g^{-1}(\overline{Y}), \widehat{\eta}_1^{(0)}(x_1) = \cdots = \widehat{\eta}_d^{(0)}(x_d) = 0$ worked well. However, a parametric fit may not be contained in the ball of Theorem 3 with probability tending to one. An alternative is to use the marginal integration estimator proposed by Linton and Härdle [11]. The latter is consistent, but costs heavier numerical calculations.

REMARK 5. If one models the conditional variance as $V(\cdot) = 1/g'(\cdot)$, then $q_2(u, y) = -[g'(g^{-1}(u))]^{-1}$. Thus, the condition for $q_2(u, y)$ is fulfilled if $g$ is strictly increasing. If one uses, as an initial estimator, the maximum smoothed quasilikelihood estimator that results from this modelling, then the *global* concavity condition on $q_2$ can be relaxed to a *local* concavity at the true function. This is because the initial estimator lies in a shrinking ball centered at the true function with probability tending to one.

4.2. *Local linear smooth backfitting.* Here, we present the theory for the maximum smoothed quasilikelihood estimator $\widehat{\boldsymbol{\eta}}$ based on local linear fit. We recall that, in the local linear case, $\widehat{\boldsymbol{\eta}}(\mathbf{x}) = (\widehat{\eta}_0(\mathbf{x}), \widehat{\eta}_1(x_1), \ldots, \widehat{\eta}_d(x_d))^T$ and $\widehat{\eta}_0(\mathbf{x}) = \widehat{\eta}_{00} + \widehat{\eta}_{01}(x_1) + \cdots + \widehat{\eta}_{0d}(x_d)$. Also, note that $\widehat{\eta}_j(x_j)$, for $1 \leq j \leq d$, estimate $\eta_j^*(x_j) = h_j\eta_{0j}^{*\prime}(x_j) = h_j(\partial\eta_{0j}^*(x_j)/\partial x_j)$.

THEOREM 5 (Rates of convergence). *Suppose that conditions* A1–A5 *hold. Then*

$$\|\widehat{\eta}_j - \eta_j^*\|_p = O_p(n^{-2/5}), \qquad\qquad j = 0, \ldots, d,$$

$$\sup_{x_j \in [2h_j, 1-2h_j]} |\widehat{\eta}_{0j}(x_j) - \eta_{0j}^*(x_j)| = O_p(n^{-2/5}\sqrt{\log n}), \qquad j = 1, \ldots, d,$$

$$\sup_{x_j \in [2h_j, 1-2h_j]} |\widehat{\eta}_j(x_j) - \eta_j^*(x_j)| = O_p(n^{-2/5}\sqrt{\log n}), \qquad j = 1, \ldots, d.$$

The asymptotic distribution of the local linear maximum smoothed quasi-likelihood estimator is given below. To state the theorem, define

$$\beta_j(x_j) = \tfrac{1}{2}\delta_j^2\left(\int t^2 K^0(t)\,dt\right)\eta_{0j}^{*\prime\prime}(x_j),$$

$$\beta_0 = -\left(\int w^*(\mathbf{x})\,d\mathbf{x}\right)^{-1}\left[\sum_{k=1}^{d}\tfrac{1}{2}\delta_k^2\int t^2 K^0(t)\,dt\int \eta_{0k}^{*\prime\prime}(x_k)w_k^*(x_k)\,dx_k\right.$$

$$\left.+\sum_{k=1}^{d}\delta_k^2\kappa(\eta_{0k}^{*\prime}(0+)w_k^*(0+)-\eta_{0k}^{*\prime}(1-)w_k^*(1-))\right].$$

Let $v_j(x_j)$ be defined as in (19).

THEOREM 6 (Asymptotic distributions). *Under the conditions of Theorem* 1, $n^{2/5}(\widehat{\eta}_{00}-\eta_{00}^*)\xrightarrow{p}\beta_0$, *and for any* $x_1,\ldots,x_d\in(0,1)$, $n^{2/5}[\widehat{\eta}_{01}(x_1)-\eta_{01}^*(x_1),\ldots,\widehat{\eta}_{0d}(x_d)-\eta_{0d}^*(x_d)]^T$ *converges in distribution to the $d$-variate Normal distribution with mean vector* $[\beta_1(x_1),\ldots,\beta_d(x_d)]^T$ *and variance matrix* $\mathrm{diag}[v_j(x_j)]$.

Theorem 6 tells that our local linear maximum smoothed quasilikelihood estimator has the oracle bias as well as the oracle variance. This property is shared with the local linear smooth backfitting estimator in the ordinary additive models. It may be interesting to compare the bias and variance properties of our estimator with those of the two-stage estimator proposed by Horowitz and Mammen [8]. Each estimator achieves the bias of the respective oracle estimator based on knowing all other components. As for the variances, we note that if the conditional density $f_{Y|\mathbf{X}}(y|\mathbf{x})$ of $Y$ given $\mathbf{X}=\mathbf{x}$ belongs to an exponential family, that is, $f_{Y|\mathbf{X}}(y|\mathbf{x})=\exp[\frac{y\theta(\mathbf{x})-b(\theta(\mathbf{x}))}{a(\phi)}+c(y,\phi)]$ for known functions $a,b$ and $c$, and one uses the canonical link $g=(b')^{-1}$, then the asymptotic bias of the two-stage estimator equals

$$v_j^{\mathrm{HM}}(x_j) = a(\phi)[E(b''(\eta^*(\mathbf{X}))^2|X_j=x_j)]^{-2}E(b''(\eta^*(\mathbf{X}))^3|X_j=x_j)$$

$$\times p_j(x_j)^{-1}\delta_j^{-1}\int K^0(t)^2\,dt.$$

An application of Hölder inequality shows that $v_j^{\mathrm{HM}}(x_j)\geq v_j(x_j)$.

THEOREM 7 (Convergence of outer and inner iterations). *Under conditions* A1–A5, *Theorems* 3 *and* 4 *remain valid for the outer and inner iterations to compute the local linear maximum smoothed quasilikelihood estimator, with* $\widehat{\boldsymbol{\eta}}$ *and* $\widehat{\boldsymbol{\eta}}^{(k)}$ *being now replaced by* $\widehat{\boldsymbol{\eta}}_L$ *and* $\widehat{\boldsymbol{\eta}}_L^{(k)}$, *respectively.*

**5. Numerical properties.** We compared our maximum smoothed likelihood estimators (YPM) with the two-stage procedures of Horowitz and Mammen [8], denoted by HM. These numerical experiments were done by R on Windows. For HM, we used R function `bs()` in the library `gam` to generate B-splines, and `nlm()` for the optimization in the first stage.

The simulation was done under the following two models for the conditional distribution:

1. $Y|\mathbf{X} = \mathbf{x} \sim \text{Bernoulli}(m(\mathbf{x}))$, where $\text{logit}(m(\mathbf{x})) = \sin(\pi x_1) + 0.5[x_2 + \sin(\pi x_2)]$;
2. $Y|\mathbf{X} = \mathbf{x} \sim \text{Poisson}(m(\mathbf{x}))$, where $\log(m(\mathbf{x})) = \sin(\pi x_1) + 0.5[x_2 + \sin(\pi x_2)]$.

We considered the following two models for the covariate vector $(X_1, X_2)$:

1. $(X_1, X_2)$ have $N_2(0, 0; 1, 1, 0)$ distribution truncated on $[-1, 1]^2$,
2. $(X_1, X_2)$ have $N_2(0, 0; 1, 1, 0.9)$ distribution truncated on $[-1, 1]^2$,

where $N_2(\mu_1, \mu_2; \sigma_1^2, \sigma_2^2, \rho)$ denotes the bivariate normal distribution with means $\mu_1, \mu_2$, variances $\sigma_1^2, \sigma_2^2$, and correlation coefficient $\rho$. Because of the truncation, the actual correlation coefficient in the second model equals 0.682. We call these models, Model $(i, j)$, where $i$ denotes the model number for the conditional distribution and $j$ is the model number for the marginal distribution of the covariate vector. For Models $(1, 1)$ and $(1, 2)$, the components $\eta_1^*$ and $\eta_2^*$ that satisfy the normalizing constraint given at (5) are $\eta_1^*(x_1) = \cos(\pi x_1)$ and $\eta_2^*(x_2) = 0.5[x_2 + \sin(\pi x_2)]$ so that $\eta_0^* = 0$. For Model $(2, 1)$, they are $\eta_1^*(x_1) = \cos(\pi x_1) - 0.4533$ and $\eta_2^*(x_2) = 0.5[x_2 + \sin(\pi x_2)] - 0.3230$ so that $\eta_0^* = 0.7763$, and for Model $(2, 2)$, they are $\eta_1^*(x_1) = \cos(\pi x_1) - 0.5874$ and $\eta_2^*(x_2) = 0.5[x_2 + \sin(\pi x_2)] - 0.4536$ so that $\eta_0^* = 1.0410$.

We generated 1,000 pseudo samples of sizes $n = 100, 500$ from each model. All the integrals involved in the smooth backfitting procedure were calculated by a trapezoidal rule based on 41 equally spaced grid points on $[-1, 1]$ for each direction. We used the theoretically optimal bandwidths for YPM. For the implementation of HM, one needs to choose the numbers of knots $\kappa_i$ at the first stage and the bandwidths at the second stage. We chose $\kappa_1 = \kappa_2 = 2$ for $n = 100$ and $\kappa_1 = \kappa_2 = 4$ for $n = 500$. We used the same bandwidths as in YPM. In a preliminary experiment with HM, we found that HM was unstable at the second stage. In our simulation, we applied a modified version of the second stage procedure, dropping the second term in the second derivative of the weighted sum of the squared errors, $S''_{nj1}(x^1, \tilde{m})$ in their notation.

Table 1 summarizes the results of the experiments. It contains the average values, over the two components, of the integrated squared biases (ISB), the integrated variances (IV) and the mean integrated squared errors (MISE), of the estimators. Note that the target components of HM are different from those of YPM by constants. This is because HM uses a different normalizing constraint that the mean of each component function is zero. The results in Table 1 are with respect to their respective targets. In calculation of the values in Table 1, we excluded bad estimates whose squared $L_2$ distance was greater than 50, that is, $\|\hat{\eta} - \eta^*\|_2^2 > 50$. In fact, HM often produced bad estimates for $n = 100$ when the covariates are correlated. Table 2 reports the number of bad estimates out of 1,000.

Table 1

*Average values, over the two components, of the integrated squared biases (ISB), the integrated variances (IV) and the mean integrated squared errors (MISE) of the maximum smoothed quasilikelihood estimator (YPM) and the two-stage estimator of Horowitz and Mammen (HM), based on 1,000 samples for the four models given in the text. LC stands for the estimators based on Nadaraya–Watson smoothing, and LL for the estimators based on local linear fit*

|  |  | $n = 100$ | | | | $n = 500$ | | | |
|---|---|---|---|---|---|---|---|---|---|
| **Model** | | **YPM LC** | **HM LC** | **YPM LL** | **HM LL** | **YPM LC** | **HM LC** | **YPM LL** | **HM LL** |
| $(1, 1)$ | ISB | 0.098 | 0.081 | 0.044 | 0.057 | 0.045 | 0.044 | 0.021 | 0.020 |
| | IV | 0.145 | 0.448 | 0.340 | 0.895 | 0.040 | 0.041 | 0.074 | 0.077 |
| | MISE | 0.243 | 0.529 | 0.384 | 0.952 | 0.084 | 0.085 | 0.095 | 0.096 |
| $(2, 1)$ | ISB | 0.068 | 0.122 | 0.023 | 0.052 | 0.017 | 0.026 | 0.009 | 0.011 |
| | IV | 0.068 | 0.371 | 0.137 | 0.545 | 0.020 | 0.020 | 0.023 | 0.023 |
| | MISE | 0.136 | 0.492 | 0.161 | 0.597 | 0.037 | 0.046 | 0.032 | 0.033 |
| $(1, 2)$ | ISB | 0.134 | 0.185 | 0.047 | 0.061 | 0.052 | 0.071 | 0.017 | 0.019 |
| | IV | 0.191 | 1.397 | 0.486 | 2.826 | 0.054 | 0.279 | 0.142 | 0.366 |
| | MISE | 0.325 | 1.581 | 0.533 | 2.887 | 0.106 | 0.349 | 0.158 | 0.385 |
| $(2, 2)$ | ISB | 0.098 | 0.170 | 0.033 | 0.143 | 0.033 | 0.041 | 0.007 | 0.014 |
| | IV | 0.125 | 1.061 | 0.370 | 2.059 | 0.027 | 0.277 | 0.054 | 0.275 |
| | MISE | 0.223 | 1.231 | 0.403 | 2.202 | 0.060 | 0.317 | 0.061 | 0.289 |

Table 2

*Number of bad estimates out of 1,000 replications for the four models given in the text ($d = 2$)*

|  | $n = 100$ | | | | $n = 500$ | | | |
|---|---|---|---|---|---|---|---|---|
| **Model** | **YPM LC** | **HM LC** | **YPM LL** | **HM LL** | **YPM LC** | **HM LC** | **YPM LL** | **HM LL** |
| $(1, 1)$ | 0 | 32 | 0 | 8 | 0 | 0 | 0 | 0 |
| $(2, 1)$ | 0 | 74 | 0 | 38 | 0 | 0 | 0 | 0 |
| $(1, 2)$ | 0 | 164 | 0 | 152 | 0 | 8 | 0 | 8 |
| $(2, 2)$ | 0 | 282 | 13 | 175 | 0 | 58 | 0 | 13 |

TABLE 3

*Average values, over the first two components, of the integrated squared biases (ISB), the integrated variances (IV) and the mean integrated squared errors (MISE), based on 100 samples of size $n = 500$ for $d = 2, 5$ and for the Bernoulli model with*
$$\text{logit}(m(\mathbf{x})) = \sin(\pi x_1) + 0.5[x_2 + \sin(\pi x_2)] + 0.1 \sum_{j=3}^{d} x_j$$

| | | (1, 1) | | | | (1, 2) | | | |
|---|---|---|---|---|---|---|---|---|---|
| | | **YPM LC** | **HM LC** | **YPM LL** | **HM LL** | **YPM LC** | **HM LC** | **YPM LL** | **HM LL** |
| $d$ | | | | | | | | | |
| 2 | ISB | 0.045 | 0.044 | 0.021 | 0.020 | 0.052 | 0.071 | 0.017 | 0.019 |
| | IV | 0.040 | 0.041 | 0.074 | 0.077 | 0.054 | 0.279 | 0.142 | 0.366 |
| | MISE | 0.084 | 0.085 | 0.095 | 0.096 | 0.106 | 0.349 | 0.158 | 0.385 |
| 5 | ISB | 0.068 | 0.041 | 0.047 | 0.014 | 0.065 | 0.179 | 0.044 | 0.019 |
| | IV | 0.035 | 0.080 | 0.093 | 0.112 | 0.043 | 0.366 | 0.171 | 0.650 |
| | MISE | 0.103 | 0.121 | 0.141 | 0.126 | 0.108 | 0.545 | 0.215 | 0.669 |

Comparing YPM and HM with the results in Table 1, we see that YPM has smaller values of MISE than HM in all cases. For correlated covariates, our simulation results suggest that IV of HM gets significantly worse whereas YPM continues to have good performance. This is mostly due to the fact that HM is unstable on the boundary of the support of the covariate vector. The good performance of YPM for correlated covariates is also in accordance with that of smooth backfitting for models with the identity link, see Nielsen and Sperlich [20]. The results also reveal that HM becomes very unstable when the sample size gets smaller. Another interesting point is that while the local linear YPM and HM certainly have less bias than their local constant versions, they have increased variance in comparison with the latter.

To see whether YPM remains competitive for higher dimensional covariates, we conducted an additional simulation with the Bernoulli model for $3 \leq d \leq 5$ where $\text{logit}(m(\mathbf{x})) = \sin(\pi x_1) + 0.5[x_2 + \sin(\pi x_2)] + 0.1 \sum_{j=3}^{d} x_j$. The covariates $X_1$ and $X_2$ were the same as in Model $(1,1)$ or $(1,2)$. The additional covariates $X_j$ for $j \geq 3$ were generated from $U(-1, 1)$ independently of other covariates. The theoretically optimal bandwidths were used for $h_1$ and $h_2$, and all other bandwidths were set to 0.2. We found that YPM continues to dominate HM for all $d$ when $X_1$ and $X_2$ are correlated. We report the results for $d = 5$ only. Table 3 shows the average values, over the first two components, of ISB, IV and MISE that are based on 100 samples of size $n = 500$.

Implementation of YPM involves multiple numerical integration so that the computational costs increase as $d$ gets high. However, one may speed up the computing time for YPM by applying a well devised Monte Carlo method for the numerical integration. If one uses an efficient numerical integration method whose grid points are as many as in one-dimensional integration, the computing time $T_d$ for $d$-dimensional covariates $(d > 3)$ equals $O(d^2) \times T_3$ since the smooth backfitting requires only two-dimensional marginal values of the weight functions. Note that, for $0.26 < \alpha < 0.3$, HM needs $d \times O(n^\alpha) \equiv \lambda$-dimensional nonlinear optimization which involves iterative inversions of $\lambda \times \lambda$ matrices. This means YPM may be as fast as, or even faster than, HM with efficient numerical integration. We do not pursue this computational issue further here since it is beyond the scope of the paper. In our current computing environments with 21 grid points in each direction, the average times (in seconds) to compute YPM and HM with a sample of size $n = 500$ for Models $(1, 1)$ and $(1, 2)$ are as reported in Table 4.

**6. Proofs and technical details.** We give only proofs of Theorems 1–4. The ideas of these proofs can be carried over to those of Theorems 5–7 for the local linear maximum smoothed quasilikelihood estimator. We note that the boundary modified kernel $K_h(u, v)$ differs from the base kernel $K_h^0(u - v)$ only when $u \in [2h, 1 - 2h]^c$ and $v \in [h, 1 - h]^c$ for $h \leq 1/2$. We will use this property repeatedly in the following proofs.

We will argue that, if a point $\bar{\boldsymbol{\eta}}$ fulfills

$$(22) \qquad \|\widehat{\mathbf{F}}(\bar{\boldsymbol{\eta}})\| = O_p(\epsilon_n) \qquad [\text{or } o_p(\epsilon_n), \text{ resp.}],$$

then $\bar{\boldsymbol{\eta}}$ also satisfies

$$(23) \qquad \|\widehat{\boldsymbol{\eta}} - \bar{\boldsymbol{\eta}}\| = O_p(\epsilon_n) \qquad [\text{or } o_p(\epsilon_n), \text{ resp.}].$$

We consider two norms: $\| \cdot \|_{w^*}$ and $\| \cdot \|_\infty$, where

$$\|\boldsymbol{\eta}\|_{w^*} = \left[ \int \left( \eta_0^2 + \sum_{j=1}^d \eta_j(x_j)^2 \right) w^*(\mathbf{x}) \, d\mathbf{x} \right]^{1/2},$$

$$\|\boldsymbol{\eta}\|_\infty = \max\{|\eta_0|, \|\eta_1\|_{\infty,1}, \ldots, \|\eta_d\|_{\infty,d}\},$$

TABLE 4
*Average computing times (in seconds) for YPM and HM with 21 grid points in each direction, for the Bernoulli model with*
$\text{logit}(m(\mathbf{x})) = \sin(\pi x_1) + 0.5[x_2 + \sin(\pi x_2)] + 0.1 \sum_{j=3}^d x_j$ *and for the sample size $n = 500$*

| | (1, 1) | | | | (1, 2) | | | |
|---|---|---|---|---|---|---|---|---|
| $d$ | YPM LC | HM LC | YPM LL | HM LL | YPM LC | HM LC | YPM LL | HM LL |
| 2 | 0.38 | 0.72 | 0.87 | 0.73 | 0.41 | 0.89 | 1.06 | 0.92 |
| 3 | 0.83 | 1.08 | 2.59 | 1.11 | 0.88 | 1.40 | 3.55 | 1.42 |
| 4 | 3.02 | 3.11 | 4.63 | 3.15 | 3.47 | 3.41 | 5.62 | 3.44 |
| 5 | 6.67 | 4.93 | 14.51 | 4.94 | 9.24 | 5.31 | 19.78 | 5.33 |

and $\|g\|_{\infty,j} = \sup_{u \in \mathcal{I}_j} |g(u)|$ for $\mathcal{I}_j = [2h_j, 1 - 2h_j]$, $j = 1, \ldots, d$.

To show that (22) implies (23), we use a version of the Newton–Kantorovich theorem. Let $\mathcal{X}$ and $\mathcal{Y}$ be Banach spaces, $F$ be a mapping $B_r(\boldsymbol{\zeta}_0) \subset \mathcal{X} \to \mathcal{Y}$, where $B_r(\boldsymbol{\zeta}_0)$ denotes a ball centered at $\boldsymbol{\zeta}_0$ with radius $r$, and $F'$ be the Fréchet derivative of $F$.

PROPOSITION 1 (Newton–Kantorovich). *Suppose that there exist constants $\alpha$, $\beta$, $c$ and $r$ such that $2\alpha\beta c < 1$ and $2\alpha < r$ for which $F$ has a derivative $F'(\boldsymbol{\zeta})$ for $\boldsymbol{\zeta} \in B_r(\boldsymbol{\zeta}_0)$, $F'$ is invertible, $\|F'(\boldsymbol{\zeta}_0)^{-1}F(\boldsymbol{\zeta}_0)\| \leq \alpha$, $\|F'(\boldsymbol{\zeta}_0)^{-1}\| \leq \beta$, $\|F'(\boldsymbol{\zeta}) - F'(\boldsymbol{\zeta}')\| \leq c\|\boldsymbol{\zeta} - \boldsymbol{\zeta}'\|$ for all $\boldsymbol{\zeta}$, $\boldsymbol{\zeta}' \in B_r(\boldsymbol{\zeta}_0)$. Then $F(\boldsymbol{\zeta}) = 0$ has a unique solution $\boldsymbol{\zeta}^*$ in $B_{2\alpha}(\boldsymbol{\zeta}_0)$. Furthermore, $\boldsymbol{\zeta}^*$ can be approximated by Newton's iterative method $\boldsymbol{\zeta}_{k+1} = \boldsymbol{\zeta}_k - F'(\boldsymbol{\zeta}_k)^{-1}F(\boldsymbol{\zeta}_k)$, which converges at a geometric rate: $\|\boldsymbol{\zeta}_k - \boldsymbol{\zeta}^*\| \leq \alpha 2^{-(k-1)}q^{2^k-1}$ where $q = 2\alpha\beta c < 1$.*

For the proof and technical details of the proposition, see Deimling [3], Section 15, for example. Proposition 1 has two important implications. One is that the distance between the unique solution and the initial point is less than $2\alpha$. This shows that (22) implies (23). The other is that, if one has a good initial guess satisfying the sufficient conditions of the proposition then one can obtain the unique solution of the equation by using the iterative method which converges geometrically fast.

We apply the proposition with $F = \widehat{\mathbf{F}}$ for the proofs of Theorems 1–3. For Theorem 1, we take $\boldsymbol{\zeta}_0 = \boldsymbol{\eta}^*$. For Theorem 2, we put $\boldsymbol{\zeta}_0$ to be some relevant approximation of $\widehat{\boldsymbol{\eta}}$. For Theorem 3, we work with $\boldsymbol{\zeta}_0 = \widehat{\boldsymbol{\eta}}^{(0)}$. For the proofs of Theorems 1 and 2, we need the following series of lemmas.

LEMMA 1. *Under conditions of Theorem 1, we have*

$$\|\widehat{\mathbf{F}}(\boldsymbol{\eta}^*)\|_{w^*} = O_p(n^{-2/5}) \quad and \quad \|\widehat{\mathbf{F}}(\boldsymbol{\eta}^*)\|_\infty = O_p(n^{-2/5}\sqrt{\log n}).$$

PROOF. Let $\psi(u) = -q_2(u, g^{-1}(u)) = [V(g^{-1}(u))g'(g^{-1}(u))]^{-1}$. With a Taylor expansion, we have for $j = 1, \ldots, d$

$$E \int \widetilde{m}(\mathbf{x})\psi(\eta^*(\mathbf{x}))\widehat{p}(\mathbf{x}) \, d\mathbf{x}_{-j}$$

$$= \int g^{-1}(\eta^*(\mathbf{x}))\psi(\eta^*(\mathbf{x}))p(\mathbf{x}) \, d\mathbf{x}_{-j} + R_{1,j,n}(x_j),$$

$$E \int g^{-1}(\eta^*(\mathbf{x}))\psi(\eta^*(\mathbf{x}))\widehat{p}(\mathbf{x}) \, d\mathbf{x}_{-j}$$

$$= \int g^{-1}(\eta^*(\mathbf{x}))\psi(\eta^*(\mathbf{x}))p(\mathbf{x}) \, d\mathbf{x}_{-j} + R_{2,j,n}(x_j),$$

where the remainders $R_{i,j,n}$ for $i = 1, 2$ and $j = 1, \ldots, d$ satisfy

$$
\begin{align}
&\sup\{|R_{i,j,n}(x_j)| : x_j \in [h_j, 1 - h_j]\} \leq (\text{const.}) n^{-2/5}, \tag{24}\\
&\sup\{|R_{i,j,n}(x_j)| : x_j \in [0, h_j) \cup (1 - h_j, 1]\} \leq (\text{const.}) n^{-1/5}.
\end{align}
$$

The above inequalities are consequences of the standard theory of kernel smoothing and properties of the boundary corrected kernels.

Since $\prod_{l \neq j}^{d} K_{h_l}(x_l, X_l^1) = \prod_{l \neq j}^{d} K_{h_l}^0(x_l - X_l^1)$ when $X_l^1 \in [h_l, 1 - h_l]$ for all $l \neq j$, and thus

$$
\left| \int g^{-1}(\eta^*(\mathbf{x})) \psi(\eta^*(\mathbf{x})) \left( \prod_{l \neq j}^{d} K_{h_l}(x_l, X_l^1) - \prod_{l \neq j}^{d} K_{h_l}^0(x_l - X_l^1) \right) d\mathbf{x}_{-j} \right|
$$

$$
\leq (\text{const.}) \sum_{l \neq j}^{d} I(X_l^1 \in [h_l, 1 - h_l]^c),
$$

we obtain

$$
\begin{align}
&\operatorname{var}\left[ \int g^{-1}(\eta^*(\mathbf{x})) \psi(\eta^*(\mathbf{x})) \widehat{p}(\mathbf{x}) \, d\mathbf{x}_{-j} \right] \notag\\
&\quad = n^{-1} \operatorname{var}\left[ \int g^{-1}(\eta^*(\mathbf{x})) \psi(\eta^*(\mathbf{x})) \prod_{l \neq j}^{d} K_{h_l}^0(x_l - X_l^1) \, d\mathbf{x}_{-j} K_{h_j}(x_j, X_j^1) \right] \tag{25}\\
&\qquad + o(n^{-1} h_j^{-1}) \notag\\
&\quad \leq (\text{const.}) n^{-1} h_j^{-1} + o(n^{-1} h_j^{-1}). \notag
\end{align}
$$

We also have

$$
\operatorname{var}\left[ Y^1 K_{h_j}(x_j, X_j^1) \int \psi(\eta^*(\mathbf{x})) \prod_{l \neq j}^{d} K_{h_l}(x_l, X_l^1) \, d\mathbf{x}_{-j} \right] = O(h_j^{-1}). \tag{26}
$$

From the inequalities (24)–(26), we obtain for $j = 1, \ldots, d$

$$
\left\| \int (\widetilde{m}(\mathbf{x}) - g^{-1}(\eta^*(\mathbf{x}))) \psi(\eta^*(\mathbf{x})) \widehat{p}(\mathbf{x}) \, d\mathbf{x}_{-j} \right\|_{w^*} = O_p(n^{-2/5}).
$$

Similarly, we find $\int [\widetilde{m}(\mathbf{x}) - g^{-1}(\eta^*(\mathbf{x}))] \psi(\eta^*(\mathbf{x})) \widehat{p}(\mathbf{x}) \, d\mathbf{x} = O_p(n^{-2/5})$. This concludes the proof of the first part.

For the proof of the second part, let

$$
A_{j,n}(x_j) = \int [\widetilde{m}_j(x_j) - g^{-1}(\eta^*(\mathbf{x}))] \psi(\eta^*(\mathbf{x})) \widehat{p}(\mathbf{x}) \, d\mathbf{x}_{-j}.
$$

Since $V$ and $g'$ are strictly positive and thus $|EA_{j,n}(x_j)| < (\text{const.})n^{-2/5}$ on $\mathcal{I}_j$, it suffices to show

$$(27) \qquad \sup_{x_j \in \mathcal{I}_j} |A_{j,n}(x_j) - E[A_{j,n}(x_j)|\mathbf{X}^1,\ldots,\mathbf{X}^n]| = O_p\left(\sqrt{\frac{\log n}{nh_j}}\right),$$

$$(28) \qquad \sup_{x_j \in \mathcal{I}_j} |E[A_{j,n}(x_j)|\mathbf{X}^1,\ldots,\mathbf{X}^n] - EA_{j,n}(x_j)| = O_p\left(\sqrt{\frac{\log n}{nh_j}}\right).$$

Define

$$D_n^i(x_j) = \int \psi(\eta^*(\mathbf{x})) \prod_{l \neq j}^d K_{h_l}(x_l, X_l^i)\, d\mathbf{x}_{-j}.$$

Then, for $x_j \in \mathcal{I}_j$

$$B_{j,n}(x_j) \equiv A_{j,n}(x_j) - E[A_{j,n}(x_j)|\mathbf{X}^1,\ldots,\mathbf{X}^n]$$

$$= n^{-1} \sum_{i=1}^n [Y^i - g^{-1}(\eta^*(\mathbf{X}^i))] D_n^i(x_j) K_{h_j}^0(x_j - X_j^i).$$

Let $\epsilon^i = Y^i - g^{-1}(\eta^*(\mathbf{X}^i))$ and $\widetilde{\epsilon}^i = \epsilon^i I(|\epsilon^i| \leq n^\alpha)$ for some $\alpha$ such that $r_0^{-1} < \alpha < 2/5$ where $r_0 > 5/2$ is the positive number in the condition A2. Let

$$\widetilde{B}_{j,n}(x_j) = \frac{1}{n} \sum_{i=1}^n [\widetilde{\epsilon}^i D_n^i(x_j) K_{h_j}^0(x_j - X_j^i) - E\widetilde{\epsilon}^i D_n^i(x_j) K_{h_j}^0(x_j - X_j^i)].$$

It is easy to see that $|E(\widetilde{\epsilon}^1 K_{h_j}^0(x_j - X_j^1) D_n^1(x_j))| < (\text{const.}) n^{-\alpha(r_0-1)} h_j^{-1} = o(n^{-2/5})$ uniformly over $x_j \in I_j$. Also, for an arbitrary positive sequence $\{a_n\}$, we have

$$P\left[\sup_{x_j \in \mathcal{I}_j} \left|\frac{1}{n} \sum_{i=1}^n (\epsilon^i - \widetilde{\epsilon}^i) K_{h_j}^0(x_j - X_j^i) D_n^i(x_j)\right| > a_n\right]$$

$$\leq P\left[\max_{1 \leq i \leq n} |\epsilon^i| > n^\alpha\right]$$

$$\leq nP[|\epsilon^1| > n^\alpha] \leq (\text{const.}) n^{-r_0\alpha+1} = o(1).$$

This implies $\sup_{x_j \in \mathcal{I}_j} |B_{j,n}(x_j) - \widetilde{B}_{j,n}(x_j)| = o_p(n^{-2/5})$.

Thus, to prove (27) it suffices to establish that

$$(29) \qquad \sup_{x_j \in \mathcal{I}_j} P\left[|\widetilde{B}_{j,n}(x_j)| > C\sqrt{\frac{\log n}{nh_j}}\right] \leq 2n^{-C+c_0}$$

for all $C > 0$ and a fixed constant $c_0$. The inequality (29) can be proved by a simple application of Markov inequality as in the proof of Theorem 6.1 in

Mammen and Park [18]. Proof of (28) is similar. This concludes the proof of the lemma. $\quad\square$

Next, we consider an approximation of $\widehat{\boldsymbol{\eta}}$. Write $\widehat{w}^* = \widehat{w}^{\eta^*}$ for simplicity. Consider the following system of linear equations for $\zeta_0, \zeta_1(\cdot), \ldots, \zeta_d(\cdot)$ which are obtained by linearly approximating the original estimating equations around $\boldsymbol{\eta}^*$:

$$
\zeta_0 = \left( \int \widehat{w}^*(\mathbf{x})\, d\mathbf{x} \right)^{-1} \int [\widetilde{m}(\mathbf{x}) - g^{-1}(\eta^*(\mathbf{x}))] \psi(\eta^*(\mathbf{x})) \widehat{p}(\mathbf{x})\, d\mathbf{x},
$$

(30)

$$
\zeta_j(x_j) = \widetilde{\zeta}_j(x_j) - \sum_{l=1, \neq j}^{d} \int \zeta_l(x_l) \frac{\widehat{w}_{jl}^*(x_j, x_l)}{\widehat{w}_j^*(x_j)}\, dx_l - \zeta_0, \qquad j = 1, \ldots, d,
$$

where

$$
\widetilde{\zeta}(\mathbf{x}) = [\widetilde{m}(\mathbf{x}) - g^{-1}(\eta^*(\mathbf{x}))] \psi(\eta^*(\mathbf{x})) \frac{\widehat{p}(\mathbf{x})}{\widehat{w}^*(\mathbf{x})},
$$

$$
\widetilde{\zeta}_j(x_j) = \left( \int \widehat{w}^*(\mathbf{x})\, d\mathbf{x}_{-j} \right)^{-1} \int \widetilde{\zeta}(\mathbf{x}) \widehat{w}^*(\mathbf{x})\, d\mathbf{x}_{-j}, \qquad j = 1, \ldots, d.
$$

Let $\widehat{\zeta}_0, \widehat{\zeta}_1(x_1), \ldots, \widehat{\zeta}_d(x_d)$ be the solution of the system (30) subject to $\int \zeta_j(x_j) \times \widehat{w}^*(\mathbf{x})\, d\mathbf{x} = 0, j = 1, \ldots, d$. Define $\widehat{\zeta}(\mathbf{x}) = \widehat{\zeta}_0 + \sum_{j=1}^{d} \widehat{\zeta}_j(x_j)$. Then $\widehat{\zeta}$ can be regarded as the minimizer in $\mathcal{H}(\widehat{w}^*)$ of

$$
\|\widetilde{\zeta} - \zeta\|_{\widehat{w}^*}^2 = \int [\widetilde{\zeta}(\mathbf{x}) - \zeta(\mathbf{x})]^2 \widehat{w}^*(\mathbf{x})\, d\mathbf{x}.
$$

For an approximation of $\widehat{\boldsymbol{\eta}}$, we take $\bar{\boldsymbol{\eta}} = \boldsymbol{\eta}^* + \widehat{\boldsymbol{\zeta}}$.

Derivation of the limiting distribution of $\bar{\boldsymbol{\eta}}$ is one of the key elements for the establishment of Theorem 2. Later, we will argue that the difference between $\bar{\boldsymbol{\eta}}$ and $\widehat{\boldsymbol{\eta}}$ is negligible by applying Proposition 1 with Lemmas 6 and 7. To derive the joint limiting distribution of $\bar{\eta}_j(x_j)$, we use the results of Mammen, Linton and Nielsen [15]. Note that a nonnegative weight function $w$ and its marginalizations $w_j$, if divided by $\int w(\mathbf{x})\, d\mathbf{x}$, can be regarded as a density function. Thus, we may have a version of Theorem 4 in Mammen, Linton and Nielsen [15] by making $\widehat{w}_{ij}^*$, $\widehat{w}_j^*$, $\widetilde{\zeta}$ and $\widetilde{\zeta}_j$, respectively, take the roles of their $\widehat{p}_{ij}$, $\widehat{p}_j$, $\widehat{m}$ and $\widehat{m}_j$. Define

$$
\alpha_{n,j}(x_j) = \left[ \frac{\partial}{\partial x_j} E(q_1(\eta^*(\mathbf{x}), g^{-1}(\eta^*(\mathbf{X}^1)))) | X_j^1 = x_j) \right] \frac{\int K_{h_j}(x_j, u)(u - x_j)\, du}{w_j(x_j) \int K_{h_j}(x_j, v)\, dv}.
$$

Put $\gamma_{n,j} \equiv 0$, $\widetilde{\zeta}_j^A(x_j) = \widetilde{\zeta}_j(x_j) - E[\widetilde{\zeta}_j(x_j) | \mathbf{X}^1, \ldots, \mathbf{X}^n]$ and $\widetilde{\zeta}_j^B(x_j) = E[\widetilde{\zeta}_j(x_j) | \mathbf{X}^1, \ldots, \mathbf{X}^n]$. We note that $\alpha_{n,j}(x_j) = 0$ for $x_j$ in the interior and equals $O(n^{-1/5})$ on the boundary. One can proceed as in the proofs of Theorems 3 and 4 of Mammen, Linton and Nielsen [15] to show the following three lemmas.

LEMMA 2. *Under the conditions of Theorem 2, the "high level conditions" of Mammen, Linton and Nielsen [15], that is, their conditions* (A1)–(A6), (A8) *and* (A9), *are satisfied with* $w_j^*, w_{ij}^*, \widehat{w}_j^*, \widehat{w}_{ij}^*, \widetilde{\zeta}, \widetilde{\zeta}_j$ *taking the roles of their* $p_j, p_{ij}, \widehat{p}_j, \widehat{p}_{ij}, \widehat{m}, \widehat{m}_j$, *respectively, and with* $\Delta_n = n^{-2/5}$, $\alpha_{n,j}(x_j)$, $\gamma_{n,j}, \widehat{\zeta}_j^A(x_j)$ *defined above and* $\beta$ *defined at* (20).

LEMMA 3. *Under the conditions of Theorem 2, it follows that for closed subsets* $S_1, \ldots, S_d$ *of* $(0,1)$

$$\sup_{x_j \in S_j} |\widehat{\zeta}_j^B(x_j) - \mu_{n,j}(x_j)| = o_p(n^{-2/5}), \qquad j = 1, \ldots, d,$$

*where* $\mu_{n,j}(x_j) = \alpha_{n,j}(x_j) + n^{-2/5}\beta_j(x_j)$.

LEMMA 4. *Under the conditions of Theorem 2, it follows that for closed subsets* $S_1, \ldots, S_d$ *of* $(0,1)$

$$\sup_{x_j \in S_j} |\widehat{\zeta}_j^A(x_j) - (\widetilde{\zeta}_j^A(x_j) - \widehat{\zeta}_0^A)| = o_p(n^{-2/5}), \qquad j = 1, \ldots, d,$$

*where* $\widehat{\zeta}_0^A = (\int \widehat{w}^*(\mathbf{x})\,d\mathbf{x})^{-1} n^{-1} \sum_{i=1}^n \int [Y^i - g^{-1}(\eta^*(\mathbf{X}^i))]\psi(\eta^*(\mathbf{x}))K_{\mathbf{h}}(\mathbf{x}, \mathbf{X}^i)\,d\mathbf{x}$.

From Lemmas 3 and 4, we obtain the asymptotic distribution of $\bar{\boldsymbol{\eta}}$ as is given in the following lemma.

LEMMA 5. *Under the conditions of Theorem 2,* $n^{2/5}(\bar{\eta}_0 - \eta_0^*) \xrightarrow{p} \beta_0$, *and for any* $x_1, \ldots, x_d \in (0,1)$, $n^{2/5}[\bar{\eta}_1(x_1) - \eta_1^*(x_1), \ldots, \bar{\eta}_d(x_d) - \eta_d^*(x_d)]^T$ *converges in distribution to the d-variate Normal distribution with mean vector* $[\beta_1(x_1), \ldots, \beta_d(x_d)]^T$ *and variance matrix* $\mathrm{diag}[v_j(x_j)]$.

LEMMA 6. *Under the conditions of Theorem 2, we have*

$$\|\widehat{\mathbf{F}}(\bar{\boldsymbol{\eta}})\|_{w^*} = o_p(n^{-2/5}) \quad \text{and} \quad \|\widehat{\mathbf{F}}(\bar{\boldsymbol{\eta}})\|_\infty = o_p(n^{-2/5}).$$

PROOF. Since $\sup_{x_j \in [0,1]} |\widehat{w}_j^*(x_j) - w_j^*(x_j)| = o_p(1)$ and $w_j^*(x_j)$ is bounded away from zero, it follows from (27) that

$$\sup_{x_j \in \mathcal{I}_j} |\widetilde{\zeta}_j^A(x_j)| = O_p(n^{-2/5}\sqrt{\log n}). \tag{31}$$

Since $\widehat{\zeta}_0^A = O_p(n^{-1/2})$, Lemma 4 and (31) imply

$$\sup_{x_j \in \mathcal{I}_j} |\widehat{\zeta}_j^A(x_j)| = O_p(n^{-2/5}\sqrt{\log n}). \tag{32}$$

Now, by a Taylor expansion and the definition of $\bar{\eta}_0$, it can be shown that

$$\int [\widetilde{m}(\mathbf{x}) - g^{-1}(\bar{\eta}(\mathbf{x}))]\psi(\bar{\eta}(\mathbf{x}))\widehat{p}(\mathbf{x})\,d\mathbf{x} = o_p(n^{-2/5}). \tag{33}$$

Also, a first-order approximation gives

$$\int [\widetilde{m}_j(x_j) - g^{-1}(\bar{\eta}(\mathbf{x}))]\psi(\bar{\eta}(\mathbf{x}))\widehat{p}(\mathbf{x})\,d\mathbf{x}_{-j}$$

$$= \widehat{w}_j^*(x_j)\widetilde{\zeta}_j(x_j) - \widehat{w}_j^*(x_j)\widehat{\zeta}_j(x_j) - \sum_{l=1,\neq j}^d \int \widehat{\zeta}_l(x_l)\widehat{w}_{jl}^*(x_j,x_l)\,dx_l$$

$$\qquad - \widehat{w}_j^*(x_j)\widehat{\zeta}_0 + r_{j,n}(x_j)$$

$$= r_{j,n}(x_j),$$

where, with $\bar{\bar{\eta}} \in (\boldsymbol{\eta}^*, \bar{\boldsymbol{\eta}})$,

$$r_{j,n}(x_j) = \int [\widehat{w}^{\eta^*}(\mathbf{x}) - \widehat{w}^{\bar{\eta}}(\mathbf{x})]\,d\mathbf{x}_{-j}\widehat{\zeta}_j(x_j)$$

$$\qquad + \sum_{l=1,\neq j}^d \int [\widehat{w}^{\eta^*}(\mathbf{x}) - \widehat{w}^{\bar{\eta}}(\mathbf{x})]\widehat{\zeta}_l(x_l)\,d\mathbf{x}_{-j}$$

$$\qquad + \int [\widehat{w}^{\eta^*}(\mathbf{x}) - \widehat{w}^{\bar{\eta}}(\mathbf{x})]\widehat{\zeta}_0.$$

With the smoothness conditions, it follows from Lemma 5 and (31) that

$$(34) \qquad \sup_{x_j \in \mathcal{I}_j} |r_{j,n}(x_j)| = o_p(n^{-2/5}).$$

Since $\|\cdot\|_{w^*} \leq (\text{const.})\|\cdot\|_{\infty}$, (33) and (34) complete the proof of the lemma. $\square$

To prove Theorems 1 and 2, it only remains to check the sufficient conditions in Proposition 1 for $\boldsymbol{\zeta}_0 = \boldsymbol{\eta}^*$ and $\bar{\boldsymbol{\eta}}$. We note that the arguments employed in Mammen and Nielsen [17] for a related problem can not be used here for general dimension $d$. Below in Lemma 7, we present a verification of the sufficient conditions that apply for general $d$.

LEMMA 7. *Under conditions* A1–A5, *the sufficient conditions in Proposition* 1 *hold with probability tending to one for* $F = \widehat{\mathbf{F}}$ *and* $\boldsymbol{\zeta}_0$ *being either* $\boldsymbol{\eta}^*$ *or* $\bar{\boldsymbol{\eta}}$, *with respect to the norms* $\|\cdot\|_{w^*}$ *and* $\|\cdot\|_{\infty}$.

PROOF. The Fréchet derivative of $\widehat{\mathbf{F}}$ at $\boldsymbol{\eta}$ is given by

$$(35) \qquad -\widehat{\mathbf{F}}'(\boldsymbol{\eta})\mathbf{g}(\mathbf{x}) = \begin{pmatrix} g_0 \int w^{\eta}(\mathbf{x})\,d\mathbf{x} \\ \int [g_0 + g_1(x_1) + \cdots + g_d(x_d)]\widehat{w}^{\eta}(\mathbf{x})\,d\mathbf{x}_{-1} \\ \vdots \qquad \vdots \qquad \vdots \\ \int [g_0 + g_1(x_1) + \cdots + g_d(x_d)]\widehat{w}^{\eta}(\mathbf{x})\,d\mathbf{x}_{-d} \end{pmatrix}.$$

Since $\widehat{p}$ converges to $p$, $\widehat{\mathbf{F}}'(\boldsymbol{\eta}^*)\mathbf{g}$ converges to $\mathbf{F}'(\boldsymbol{\eta}^*)\mathbf{g}$ where $\mathbf{F}'(\boldsymbol{\eta})$ is defined in the same way as $\widehat{\mathbf{F}}'(\boldsymbol{\eta})$ with $\widehat{p}$ in the latter being replaced by $p$, and thus $\widehat{w}^\eta$ being replaced by $w^\eta$. Furthermore, as in the proof of Lemma 6, one can show $\widehat{\mathbf{F}}'(\bar{\boldsymbol{\eta}})\mathbf{g}$ also converges to $\mathbf{F}'(\boldsymbol{\eta}^*)\mathbf{g}$. Here, the convergence means convergence with respect to the $\|\cdot\|_{w^*}$ or $\|\cdot\|_\infty$ norm in probability.

Note that $\Pi_j^*(\cdot) \equiv [w_j^*(x_j)]^{-1} \int \cdot w^*(\mathbf{x}) \, d\mathbf{x}_{-j}$ is the projection operator onto $\mathcal{H}_j^0(w^*)$. Let $\mathbf{w}(\mathbf{x}) = (w_1^*(x_1), \ldots, w_d^*(x_d))^T$, and define

$$A = \begin{pmatrix} 1 & \mathbf{0}^T \\ \mathbf{0} & B \end{pmatrix}, \qquad B = \begin{pmatrix} \Pi_1^* & \cdots & \Pi_1^* \\ \vdots & \ddots & \vdots \\ \Pi_d^* & \cdots & \Pi_d^* \end{pmatrix},$$

$$D = \begin{pmatrix} \int w^*(\mathbf{x}) \, d\mathbf{x} & \mathbf{0}^T \\ \mathbf{w}(\cdot) & \mathrm{diag}(\mathbf{w}(\cdot)) \end{pmatrix}.$$

Then, one can write

$$\mathbf{F}'(\boldsymbol{\eta}^*)\mathbf{g} = DA\mathbf{g}.$$

We note that $D^{-1}$ is bounded since $g'(m(\cdot))^2 V(m(\cdot))$ is bounded, and $p$ is bounded away from zero.

We only need to show that the linear operator $A$ has a bounded inverse and the Lipschitz condition is satisfied for $\mathbf{F}'$. Note that the linear operator $A$ has a bounded inverse if $B$ has. We apply the inverse mapping theorem to show the linear operator $B$ has a bounded inverse. In the proof, the spaces are redefined by dropping the constant.

Suppose that $B\mathbf{g} = \mathbf{0}$ for a given $\mathbf{g} \in \mathcal{G}^0(w^*)$. Then we have $\Pi_j^*(g_1 + \cdots + g_d) = 0$ for $j = 1, \ldots, d$. This implies that

$$g_1 + \cdots + g_d \in \mathcal{H}_1^{0\perp} \cap \cdots \cap \mathcal{H}_d^{0\perp} = (\mathcal{H}_1^0 + \cdots + \mathcal{H}_d^0)^\perp.$$

Thus, $g_1 + \cdots + g_d = 0$ so that $\mathbf{g} = \mathbf{0}$. Hence, $B$ is one-to-one. Next, note that $B$ is self-adjoint, that is, $B^* = B$ since, for any $\mathbf{g}, \boldsymbol{\gamma} \in \mathcal{G}^0(w^*)$,

$$\begin{aligned}
\langle B\mathbf{g}, \boldsymbol{\gamma} \rangle &= \langle \Pi_1^*(g_1 + \cdots + g_d), \gamma_1 \rangle_{\mathcal{H}_1^0} + \cdots + \langle \Pi_d^*(g_1 + \cdots + g_d), \gamma_d \rangle_{\mathcal{H}_d^0} \\
&= \langle g_1 + \cdots + g_d, \gamma_1 \rangle_{\mathcal{H}^0} + \cdots + \langle g_1 + \cdots + g_d, \gamma_d \rangle_{\mathcal{H}^0} \\
&= \langle g_1, \gamma_1 + \cdots + \gamma_d \rangle_{\mathcal{H}^0} + \cdots + \langle g_d, \gamma_1 + \cdots + \gamma_d \rangle_{\mathcal{H}^0} \\
&= \langle \Pi_1^* g_1, \gamma_1 + \cdots + \gamma_d \rangle_{\mathcal{H}^0} + \cdots + \langle \Pi_d^* g_d, \gamma_1 + \cdots + \gamma_d \rangle_{\mathcal{H}^0} \\
&= \langle g_1, \Pi_1^*(\gamma_1 + \cdots + \gamma_d) \rangle_{\mathcal{H}_1^0} + \cdots + \langle g_d, \Pi_d^*(\gamma_1 + \cdots + \gamma_d) \rangle_{\mathcal{H}_d^0} \\
&= \langle \mathbf{g}, B\boldsymbol{\gamma} \rangle,
\end{aligned}$$

where we use the subscripts to emphasize which inner product is used. Thus, $\mathbf{R}(B)^\perp = \mathbf{N}(B^*) = \mathbf{N}(B) = \{\mathbf{0}\}$, where $\mathbf{R}$ and $\mathbf{N}$ denote the range- and null-spaces, respectively. This implies $B$ is onto.

We can conclude that $B$ has a bounded inverse if we prove $B$ is bounded. The boundedness in $\|\cdot\|_{w^*}$ of $B$ can be easily checked as follows:

$$\|B\mathbf{g}\|_{w^*}^2 = \int [\{\Pi_1^*(g_1 + \cdots + g_d)\}^2 + \cdots + \{\Pi_d^*(g_1 + \cdots + g_d)\}^2] w^*(\mathbf{x})\,d\mathbf{x}$$

$$\leq \sum_{i=1}^d \|\Pi_i^*\|_{w^*}^2 \|g_1 + \cdots + g_d\|_{w^*}^2$$

$$\leq d^2 \int [g_1(x_1)^2 + \cdots + g_d(x_d)^2] w^*(\mathbf{x})\,d\mathbf{x} = d^2 \|\mathbf{g}\|_{w^*}^2.$$

Now, for the norm $\|\cdot\|_\infty^*$ defined by

$$\|\mathbf{g}\|_\infty^* \equiv \max\left\{ \sup_{x_1 \in (0,1)} |g_1(x_1)|, \ldots, \sup_{x_d \in (0,1)} |g_d(x_d)| \right\},$$

we have $\|B\mathbf{g}\|_\infty^* \leq d\|\mathbf{g}\|_\infty^*$ since

$$\left| \int g_k(x_k) w^*(\mathbf{x})\,d\mathbf{x}_{-j} / w_j^*(x_j) \right| \leq \sup_{x_k \in (0,1)} |g_k(x_k)|.$$

To check the Lipschitz condition, we write $\mathbf{F}'(\boldsymbol{\eta}) = D(\boldsymbol{\eta})A(\boldsymbol{\eta})$ where $D(\boldsymbol{\eta})$ and $A(\boldsymbol{\eta})$ are defined in the same way as $D$ and $A$, respectively, with $\boldsymbol{\eta}$ substituting for $\boldsymbol{\eta}^*$ thus with $w^\eta$ substituting for $w^*$. Then,

$$\|\mathbf{F}'(\boldsymbol{\eta}) - \mathbf{F}'(\boldsymbol{\eta}')\|_{w^*} \leq \|D(\boldsymbol{\eta})\|_{w^*} \|A(\boldsymbol{\eta}) - A(\boldsymbol{\eta}')\|_{w^*}$$
$$+ \|D(\boldsymbol{\eta}) - D(\boldsymbol{\eta}')\|_{w^*} \|A(\boldsymbol{\eta}')\|_{w^*}.$$

Since $g'(m(\cdot))^2 V(m(\cdot))$ and $p$ are bounded away from zero and infinity, $\|D(\boldsymbol{\eta})\|_{w^*}$ and $\|A(\boldsymbol{\eta}')\|_{w^*}$ are bounded by some constant. From the smoothness of $g'(m(\cdot))^2 V(m(\cdot))$, we also have $\|(D(\boldsymbol{\eta}) - D(\boldsymbol{\eta}'))\mathbf{g}\|_{w^*} \leq$ (const.)$\|\mathbf{g}\|_{w^*}\|\boldsymbol{\eta} - \boldsymbol{\eta}'\|_{w^*}$ and $\|(A(\boldsymbol{\eta}) - A(\boldsymbol{\eta}'))\mathbf{g}\|_{w^*} \leq$ (const.)$\|\mathbf{g}\|_{w^*}\|\boldsymbol{\eta} - \boldsymbol{\eta}'\|_{w^*}$. This establishes $\|\mathbf{F}'(\boldsymbol{\eta}) - \mathbf{F}'(\boldsymbol{\eta}')\|_{w^*} \leq$ (const.)$\|\boldsymbol{\eta} - \boldsymbol{\eta}'\|_{w^*}$. Checking the Lipschitz condition for the norm $\|\cdot\|_\infty$ is similar, hence omitted. $\square$

The following lemma tells that the norms $\|\cdot\|_{w^\mu}$ and $\|\cdot\|_{\widehat{w}^\mu}$ are equivalent to $\|\cdot\|_p$ and $\|\cdot\|_{\widehat{p}}$, respectively.

LEMMA 8. *Suppose that conditions* A1–A5 *hold. For any continuous function $\mu$, there exist positive constants $c$ and $C$ such that for each $\eta \in \mathcal{H}(p)$ and $\boldsymbol{\eta} \in \mathcal{G}(p)$,*

$$c\|\eta\|_p \leq \|\eta\|_{w^\mu} \leq C\|\eta\|_p \quad and \quad c\|\boldsymbol{\eta}\|_p \leq \|\boldsymbol{\eta}\|_{w^\mu} \leq C\|\boldsymbol{\eta}\|_p.$$

*Also, there exist positive constants $c'$ and $C'$ such that for each $\eta \in \mathcal{H}(\widehat{p})$ and $\boldsymbol{\eta} \in \mathcal{G}(\widehat{p})$,*

$$c'\|\eta\|_{\widehat{p}} \leq \|\eta\|_{\widehat{w}^\mu} \leq C'\|\eta\|_{\widehat{p}} \quad and \quad c'\|\boldsymbol{\eta}\|_{\widehat{p}} \leq \|\boldsymbol{\eta}\|_{\widehat{w}^\mu} \leq C'\|\boldsymbol{\eta}\|_{\widehat{p}}$$

*with probability tending to one.*

PROOF.   From the condition A2 and the continuity of $\mu$ and $m$, the function $-q_2(\mu(\cdot), m(\cdot))$ is bounded away from zero and infinity on any compact set. Thus, there exist positive constants $c$ and $C$ such that $cp(\mathbf{x}) \leq w^\mu(\mathbf{x}) \leq Cp(\mathbf{x})$ for all $\mathbf{x} \in [0,1]^d$. This establishes the first part of the lemma. Since $\|\eta\|_{\widehat{w}^\mu}$ and $\|\eta\|_{\widehat{p}}$ converge in probability to $\|\eta\|_{w^\mu}$ and $\|\eta\|_p$, respectively, the second part of the lemma follows.  $\square$

PROOF OF THEOREM 1.   The theorem follows directly from Lemma 1, Lemma 8 and Proposition 1 with application of Lemma 7.  $\square$

PROOF OF THEOREM 2.   The theorem follows directly from Lemma 5, Lemma 6 and Proposition 1 with application of Lemma 7.  $\square$

PROOF OF THEOREM 3.   One may show, by a parallel argument with the proof of Lemma 7 substituting $\widehat{w}^{(0)}$ for $w^*$, that there exists a constant $c_0$ such that $\|\widehat{\mathbf{F}}'(\widehat{\boldsymbol{\eta}}^{(0)})^{-1}\|_{\widehat{w}^{(0)}} \leq c_0$ with probability tending to one. Since $\|\mathbf{g}\|_{\widehat{p}}$ converges in probability to $\|\mathbf{g}\|_p$ for $\mathbf{g} \in L_2(p)$, it follows from Lemma 8 that $\|\widehat{\mathbf{F}}'(\widehat{\boldsymbol{\eta}}^{(0)})^{-1}\|_p \leq c_1$ with probability tending to one for some constant $c_1$. To check the Lipschitz condition in Proposition 1 for $\widehat{\mathbf{F}}'$, one may follow the approach in the proof of Lemma 7 with the representation $\widehat{\mathbf{F}}'(\boldsymbol{\eta}) = \widehat{D}(\boldsymbol{\eta})\widehat{A}(\boldsymbol{\eta})$, where $\widehat{D}(\boldsymbol{\eta})$ and $\widehat{A}(\boldsymbol{\eta})$ are defined in the same way as $D$ and $A$, respectively, with $\widehat{w}^\eta$ substituting for $w^*$. One can prove that there exists a constant $c_2$ such that $\|\widehat{\mathbf{F}}'(\boldsymbol{\eta}) - \widehat{\mathbf{F}}'(\boldsymbol{\eta}')\|_p \leq c_2\|\boldsymbol{\eta} - \boldsymbol{\eta}'\|_p$ with probability tending to one.

Now let $\mathbf{F}$ be defined as $\widehat{\mathbf{F}}$ with $\overline{Y}, \widetilde{m}_j, \widehat{p}_j, \widehat{p}$ being replaced by $EY, m_j, p_j, p$, respectively. Then, $\widehat{\mathbf{F}}\boldsymbol{\eta}$ converges to $\mathbf{F}\boldsymbol{\eta}$ with respect to the $\|\cdot\|_p$ norm in probability, uniformly for $\boldsymbol{\eta}$ in any compact set. From this convergence, the uniform continuity of $\mathbf{F}$ and the fact $\widehat{\mathbf{F}}\widehat{\boldsymbol{\eta}} = \mathbf{0}$, it follows that there exists a positive constant $r$ such that

$$\sup_{\boldsymbol{\eta} \in B_r(\widehat{\boldsymbol{\eta}})} \|\widehat{\mathbf{F}}\boldsymbol{\eta}\|_p < \frac{1}{2c_1^2 c_2}$$

with probability tending to one, where $B_r(\widehat{\boldsymbol{\eta}})$ is a ball in $L_2(p)$. This proves that, if $\widehat{\boldsymbol{\eta}}^{(0)} \in B_r(\widehat{\boldsymbol{\eta}})$ with probability tending to one, then

$$\|\widehat{\mathbf{F}}'(\widehat{\boldsymbol{\eta}}^{(0)})^{-1}\widehat{\mathbf{F}}\widehat{\boldsymbol{\eta}}^{(0)}\|_p \leq \frac{1}{2c_1 c_2}$$

with probability tending to one. The theorem now follows from Proposition 1.  $\square$

PROOF OF THEOREM 4.   Note that $\widehat{\Pi}_j^{(k)}(\cdot) \equiv [\widehat{w}_j^{(k)}(x_j)]^{-1} \int \cdot \widehat{w}^{(k)}(\mathbf{x}) \, d\mathbf{x}_{-j}$ are Hilbert–Schmidt operators in $L_2(\widehat{w}^{(k)})$. This implies that for each $k$ there exists a stochastic $0 < \widehat{\rho}_k < 1$ such that

$$\|\widehat{\xi}^{(k),[r]} - \widehat{\xi}^{(k)}\|_{\widehat{w}^{(k-1)}} \leq \widetilde{\rho}_k^r \|\widehat{\xi}^{(k)}\|_{\widehat{w}^{(k-1)}},$$

where $\widehat{\rho}_k < 1$. See Theorem 4.B in Appendix 4 of Bickel et al. [1] for details. This establishes the first part of the theorem.

If the initial estimator $\widehat{\boldsymbol{\eta}}^{(0)}$ belongs to the ball $B_r(\widehat{\boldsymbol{\eta}})$ with probability tending to one, then Theorem 3 tells us that, with probability tending to one, $\widehat{w}^{(k)}$ converges to $\widehat{w}^{(\infty)}$ and $\widehat{\xi}^{(k)}$ to zero (in $\|\cdot\|_p$ or $\|\cdot\|_{\widehat{p}}$ norm) as $k$ goes to infinity, where $\widehat{w}^{(\infty)}$ is defined as $\widehat{w}^{(k)}$ with $\widehat{\eta}^{(k)}$ being replaced by $\widehat{\eta}$. Define $\widehat{\rho}_\infty$ as $\rho_k$ with $\widehat{w}^{(\infty)}$ substituting for $\widehat{w}^{(k-1)}$. Note that $0 < \widehat{\rho}_\infty < 1$ since $\widehat{\Pi}_j^{(\infty)}(\cdot) \equiv [\widehat{w}_j^{(\infty)}(x_j)]^{-1} \int \cdot \widehat{w}^{(\infty)}(\mathbf{x}) \, d\mathbf{x}_{-j}$ are also Hilbert–Schmidt operators in $L_2(\widehat{w}^{(\infty)})$. This implies that within an event of probability tending to one, there exists $0 < \widehat{\rho} < 1$ and $\widehat{\varepsilon} > 0$ such that $\widehat{\rho}_k \leq \widehat{\rho}$ and $\|\widehat{\xi}^{(k)}\|_{\widehat{p}} \leq \widehat{\varepsilon}$ for all $k$. Thus, from Lemma 8 we conclude that with probability tending to one there exist $0 < \widehat{\rho} < 1$ and $\widehat{C}$, which are independent of $k$, such that

$$\|\widehat{\xi}^{(k),[r]} - \widehat{\xi}^{(k)}\|_{\widehat{p}} \leq \widehat{C}\widehat{\rho}^r.$$

This completes the proof of Theorem 4.  □

**Acknowledgments.**  We thank the two referees and an associate editor for their helpful comments on the earlier version of the paper.

## REFERENCES

[1] BICKEL, P., KLAASSEN, A., RITOV, Y. and WELLNER, J. (1993). *Efficient and Adaptive Estimation for Semiparametric Models.* The Johns Hopkins Univ. Press, Baltimore. MR1245941

[2] BUJA, A., HASTIE, T. and TIBSHIRANI, R. (1989). Linear smoothers and additive models (with discussion). *Ann. Statist.* **17** 453–510. MR0994249

[3] DEIMLING, K. (1985). *Nonlinear Functional Analysis.* Springer, Berlin. MR0787404

[4] EILERS, P. H. C. and MARX, B. D. (2002). Generalized linear additive smooth structures. *J. Comput. Graph. Statist.* **11** 758–783. MR1944262

[5] FAN, J., HECKMAN, N. E. and WAND, M. P. (1995). Local polynomial kernel regression for generalized linear models and quasi-likelihood functions. *J. Amer. Statist. Assoc.* **90** 141–150. MR1325121

[6] FRIEDMAN, J. and STUETZLE, W. (1981). Projection pursuit regression. *J. Amer. Statist. Assoc.* **76** 376, 817–823. MR0650892

[7] HASTIE, T. J. and TIBSHIRANI, R. J. (1990). *Generalized Additive Models.* Chapman and Hall, London. MR1082147

[8] HOROWITZ, J. and MAMMEN, E. (2004). Nonparametric estimation of an additive model with a link function. *Ann. Statist.* **32** 2412–2443. MR2153990

[9] KAUERMANN, G. and OPSOMER, J. D. (2003). Local likelihood estimation in generalized additive models. *Scand. J. Statist.* **30** 317–337. MR1983128

[10] LEE, Y. K. (2004). On marginal integration method in nonparametric regression. *J. Korean Statist. Soc.* **33** 435–448. MR2126371

[11] LINTON, O. and HÄRDLE, W. (1996). Estimating additive regression models with known links. *Biometrika* **83** 529–540. MR1423873

[12] LINTON, O. and NIELSEN, J. P. (1995). A kernel method of estimating structured nonparametric regression based on marginal integration. *Biometrika* **82** 93–100. MR1332841

[13] LINTON, O. (2000). Efficient estimation of generalized additive nonparametric regression models. *Econometric Theory* **16** 502–523. MR1790289

[14] LUENBERGER, D. G. (1969). *Optimization by Vector Space Methods*. Wiley, New York. MR0238472

[15] MAMMEN, E., LINTON, O. and NIELSEN, J. P. (1999). The existence and asymptotic properties of a backfitting projection algorithm under weak conditions. *Ann. Statist.* **27** 1443–1490. MR1742496

[16] MAMMEN, E., MARRON, J. S., TURLACH, B. A. and WAND, M. P. (2001). A general projection framework for constrained smoothing. *Statist. Sci.* **16** 232–248. MR1874153

[17] MAMMEN, E. and NIELSEN, J. P. (2003). Generalised structured models. *Biometrika* **90** 551–566. MR2006834

[18] MAMMEN, E. and PARK, B. U. (2005). Bandwidth selection for smooth backfitting in additive models. *Ann. Statist.* **33** 1260–1294. MR2195635

[19] MAMMEN, E. and PARK, B. U. (2006). A simple smooth backfitting method for additive models. *Ann. Statist.* **34** 2252–2271. MR2291499

[20] NIELSEN, J. and SPERLICH, S. (2005). Smooth backfitting in practice. *J. Roy. Statist. Soc. Ser. B* **67** 43–61. MR2136638

[21] OPSOMER, J. D. (2000). Asymptotic properties of backfitting estimators. *J. Multivariate Anal.* **73** 166–179. MR1763322

[22] OPSOMER, J. D. and RUPPERT, D. (1997). Fitting a bivariate additive model by local polynomial regression. *Ann. Statist.* **25** 186–211. MR1429922

[23] STONE, C. J. (1985). Additive regression and other nonparametric models. *Ann. Statist.* **13** 689–705. MR0790566

[24] STONE, C. J. (1986). The dimensionality reduction principle for generalized additive models. *Ann. Statist.* **14** 590–606. MR0840516

B. U. PARK
STATISTICS SEOUL NATIONAL UNIVERSITY
SEOUL 151-747
KOREA
E-MAIL: bupark@stats.snu.ac.kr

E. MAMMEN
K. YU
DEPARTMENT OF ECONOMICS
UNIVERSITY OF MANNHEIM
L7, 3-5
688131 MANNHEIM
GERMANY
E-MAIL: emammen@rumms-uni-mannheim.de
      yukyusan@rumms-uni-mannheim.de